\documentclass[
final
]{dmtcs-episciences}

\usepackage[utf8]{inputenc}
\usepackage{subfigure}

\author{J\o rgen Bang-Jensen\affiliationmark{1}
  \and Jonas Costa Ferreira da Silva\affiliationmark{2,3}
  \and Fr\'ed\'eric Havet\affiliationmark{3}}

\title[On the inversion number of oriented graphs.]{On the inversion number of oriented graphs. \footnote{A preliminary version of this paper was published in the proceedings of the 1st International Conference on Algebras, Graphs and Ordered Sets (ALGOS 2020). }}

\affiliation{
  Department of Mathematics and Computer  Science, University of Southern Denmark, Odense, Denmark\\
Department of Mathematics, Universidade Federal do Cear\'a, Fortaleza, Brazil\\
Universit\'e C\^ote d'Azur, CNRS, Inria, I3S, Sophia Antipolis, France}
\received{2021-05-11}
\revised{2022-02-07, 2022-07-16}
\accepted{2022-10-30}

\usepackage{times,amsmath,amsthm,amssymb,graphicx,xspace,epsfig,xcolor}
\usepackage{tikz}
\usepackage{color}

\newtheorem{theorem}{Theorem}[section]
\newtheorem{corollary}[theorem]{Corollary}
\newtheorem{proposition}[theorem]{Proposition}
\newtheorem{lemma}[theorem]{Lemma}

\theoremstyle{definition}
\newtheorem{remark}[theorem]{Remark}

\newtheorem{problem}[theorem]{Problem}
\newtheorem{conjecture}[theorem]{Conjecture}

\newcommand{\ra}{\rightarrow}
\newcommand{\la}{\leftarrow}
\newcommand{\rad}{\Rightarrow}

\newcommand{\fact}[1]{{\color{blue}#1}}
\newcommand{\sfact}[1]{{\color{red}#1}}
\newcommand{\asm}[1]{{\color{orange}#1}}

\newenvironment{subproof}{\par\noindent {\it Proof}.\ }{\hfill$\lozenge$\par\vspace{11pt}}

\DeclareMathOperator{\Inv}{Inv}
\DeclareMathOperator{\inv}{inv}

\begin{document}
\publicationdetails{23}{2022}{2}{8}{7474}
\maketitle
\begin{abstract}
 Let $D$ be an oriented graph.
The {\bf inversion} of a set $X$ of vertices in $D$ consists in reversing the direction of all arcs with both ends in $X$.  The {\bf inversion number} of $D$, denoted by $\inv(D)$, is the minimum number of inversions needed to make $D$ acyclic. Denoting by   $\tau(D)$, $\tau' (D)$, and $\nu(D)$ the cycle transversal number, the cycle arc-transversal number and the cycle packing number of $D$ respectively, one shows that $\inv(D) \leq  \tau' (D)$,  $\inv(D) \leq 2\tau(D)$ and there exists a function $g$ such that $\inv(D)\leq g(\nu(D))$.
We conjecture that for any two oriented graphs $L$ and $R$, $\inv(L\ra R) =\inv(L) +\inv(R)$ where $L\ra R$ is the dijoin of $L$ and $R$. This would imply that the first two inequalities are tight. We prove this conjecture when
$\inv(L)\leq 1$ and $\inv(R)\leq 2$ and when $\inv(L) =\inv(R)=2$ and $L$ and $R$ are strongly connected.
We also show that the function $g$ of the third inequality satisfies $g(1)\leq 4$.

We then consider the complexity of deciding whether $\inv(D)\leq k$ for a given oriented graph $D$.
We show that it is NP-complete for $k=1$, which together with the above conjecture would imply that it is NP-complete for every $k$. This contrasts with a result of  Belkhechine et al. which states that deciding whether $\inv(T)\leq k$ for a given tournament $T$ is polynomial-time solvable.
\end{abstract}

\section{Introduction}

Notation not given below is consistent with \cite{bang2009}. We denote by $[k]$ the set $\{1,2, \dots, k\}$.

Making a digraph acyclic by either removing a minimum cardinality set of arcs or vertices are  important and heavily studied problems, known under the names {\sc Cycle Arc Transversal} or {\sc Feedback Arc Set}  and  {\sc Cycle Transversal} or {\sc Feedback Vertex Set}. 
A {\bf cycle transversal} or {\bf feedback vertex set} (resp.  {\bf cycle arc-transversal} or {\bf feedback arc set})  in a digraph is a set of vertices (resp. arcs) whose deletion results in an acyclic digraph. The {\bf cycle transversal number} (resp. {\bf  cycle arc-transversal number}) is the minimum size of a cycle transversal (resp. cycle arc-transversal) of $D$ and is denoted by $\tau(D)$ (resp. $\tau'(D)$). 
 It is well-known that a digraph is acyclic if and only if it admits an {\bf acyclic ordering}, that is an ordering $(v_1, \dots , v_n)$ of its vertices such that there is no {\bf backward} arc (i.e. an arc $v_jv_i$ with $i < j$). It follows that a minimum cycle arc-transversal $F$ in a digraph $D$ consists only of backward arcs with respect to any acyclic ordering of $D\setminus F$.
Thus the digraph $D'$ obtained from $D$ by  reversing the arcs of $F$ is also acyclic.
Conversely, if the digraph $D'$ obtained from $D$ by  reversing the arcs of $F$ is acyclic, then $D\setminus F$ is also trivially acyclic.
Therefore the cycle arc-transversal number of a digraph is also the minimum size of a set of arcs whose reversal makes the digraph acyclic.

It is well-known and easy to show that $\tau(D)\leq  \tau' (D)$ (just take one end-vertex of each arc in a minimum cycle arc-transversal).

Computing $\tau(D)$ and  $\tau'(D)$ are two of the first problems shown to be NP-hard listed by Karp in~\cite{karp1972}.
They also remain NP-complete in tournaments as shown by Bang-Jensen and Thomassen~\cite{bangSJDM5} and Speckenmeyer \cite{speckenmeyerLNCS411} for $\tau$, and by Alon~\cite{alonSJDM20} and Charbit, Thomass\'e, and Yeo~\cite{charbitCPC16} for $\tau'$.

 In this paper, we consider another operation, called {\bf inversion}, where we reverse all arcs of an induced subdigraph.
Let $D$ be a digraph.
The {\bf inversion} of a set $X$ of vertices consists in reversing the direction of all arcs of $D\langle X\rangle$.
We say that we {\bf invert} $X$ in $D$. The resulting digraph is denoted by $\Inv(D;X)$.  
If $(X_i)_{i\in I}$  is a family of subsets of $V(D)$, then $\Inv(D; (X_i)_{i\in I})$ is the digraph obtained after inverting the
$X_i$ one after another. Observe that this is independent of the order in which we invert the $X_i$~: $\Inv(D; (X_i)_{i\in I})$ is obtained from $D$ by reversing the arcs such that an odd number of the $X_i$ contain its two end-vertices.

Since an inversion preserves the directed cycles of length $2$, a digraph can be made acyclic only if it has no directed cycle of length 2, that is if it is an {\bf oriented graph}. Reciprocally, observe that in an oriented graph, reversing an arc $a=uv$ is the same as inverting $X_a=\{u,v\}$. Hence if $F$ is a minimum cycle arc-transversal of $D$, then $\Inv(D; (X_a)_{a\in F})$ is acyclic.

A {\bf decycling family} of an oriented graph $D$ is a family of subsets $(X_i)_{i\in I}$ of subsets of $V(D)$ such that
 $\Inv(D; (X_i)_{i\in I})$ is acyclic.
 The {\bf inversion number} of an oriented graph $D$, denoted by $\inv(D)$,  is the minimum number of inversions needed to transform $D$ into an acyclic digraph, that is, the minimum cardinality of a decycling family.
By convention, the empty digraph (no vertices) is acyclic and so has inversion number $0$.

 \medskip

\subsection{Inversion versus cycle (arc-) transversal and cycle packing}\label{subsec:versus}

One can easily obtain the following upper bounds on the inversion number in terms of the cycle transversal number and the cycle arc-transversal number. See Section~\ref{sec:easy}.

\begin{theorem}\label{thm:bound-fvs}
$\inv(D) \leq  \tau' (D)$ and  $\inv(D) \leq 2\tau(D)$ for all oriented graph $D$.
 \end{theorem}

A natural question is to ask whether these bounds are tight or not.

We denote by $\vec{C_3}$ the directed cycle of length $3$ and by $TT_n$ the transitive tournament of order $n$. The vertices of $TT_n$ are $v_1, \dots, v_n$ and its arcs $\{v_iv_j \mid i < j\}$.
The {\bf lexicographic product}  of a digraph $D$ by a digraph $H$ is the digraph $D[H]$ with vertex set $V(D)\times V(H)$
and arc set
$A(D[H])  =  \{(a,x)(b,y) \mid ab \in A(D), \mbox{\ or\ } a=b \mbox{\ and\ } xy\in A(H)\}$.
It can be seen as blowing up each vertex of $D$ by a copy of $H$.
Using boolean dimension,  Pouzet et al.~\cite{PST} proved  the following.
\begin{theorem}[Pouzet et al.~\cite{PST}]\label{thm:TT[C3]}
 $\inv(TT_n[\vec{C_3}]) = n$. 
\end{theorem}
Since $\tau'(TT_n[\vec{C_3}])=n$, this shows that the inequality $\inv(D) \leq  \tau' (D)$  of Theorem~\ref{thm:bound-fvs} is tight.

\medskip

Pouzet asked for an elementary proof of Theorem~\ref{thm:TT[C3]}.
Let $L$ and $R$ be two oriented graphs. The {\bf dijoin} from $L$ to $R$ is the oriented graph, denoted by $L\ra R$, obtained from the disjoint union of $L$ and $R$ by adding all arcs from $L$ to $R$.
Observe that $TT_n[\vec{C_3}]=\vec{C_3}\ra TT_{n-1}[\vec{C_3}]$.
So an elementary way to prove Theorem~\ref{thm:TT[C3]} would be to prove that $\inv(\vec{C_3}\ra T) = \inv(T)+1$ for all tournament $T$.

First inverting $\inv(L)$ subsets of $V(L)$ to make $L$ acyclic and then inverting $\inv(R)$ subsets of $V(R)$ to make $R$ acyclic, makes $L\ra R$ acyclic. Therefore we have the following inequality.
\begin{proposition}\label{prop:D1>D2}
 $\inv(L\ra R) \leq \inv(L) + \inv(R)$.
 \end{proposition}
 In fact, we believe that equality always holds.

\smallskip

\begin{conjecture}\label{conj:dijoin}
For any two oriented graphs, $L$ and $R$, $\inv(L\ra R) =\inv(L) + \inv(R)$.
\end{conjecture}

As observed in Proposition~\ref{prop:equiv}, this conjecture is equivalent to its restriction to tournaments.
If $\inv(L) =0$ (resp. $\inv(R)=0$), then Conjecture~\ref{conj:dijoin} holds has any decycling family of $R$ (resp. $L$) is also
a decycling family of $L\ra R$.
In Section~\ref{sec:dijoin}, we prove Conjecture~\ref{conj:dijoin} when $\inv(L)=1$ and $\inv(R)\in \{1,2\}$.
We also prove it when $\inv(L)=\inv(R)=2$ and both $L$ and $R$ are strongly connected.

\medskip

Let us now consider the inequality  $\inv(D) \leq 2\tau(D)$ of Theorem~\ref{thm:bound-fvs}.
One can see that is tight for $\tau(D)=1$. Indeed, let $V_n$ be the tournament obtained from a $TT_{n-1}$ by adding a vertex $x$ such that
$N^+(x) = \{v_i \mid i~\mbox{is odd}\}$ (and so $N^-(x) = \{v_i \mid i~\mbox{is even}\}$.
Clearly, $\tau(V_n)=1$ because $V_n-x$ is acyclic, and one can easily check that $\inv(V_n) \geq 2$ for $n\geq 5$. 
Observe that $V_5$ is strong, so by the above results, we have $\inv(V_5\ra V_5) = 4$ while $\tau(V_5\ra V_5)=2$, so the inequality  $\inv(D) \leq 2\tau(D)$ is also tight for $\tau(D)=2$.
More generally,  Conjecture~\ref{conj:dijoin} would imply that  $\inv(TT_n[V_5]) = 2n$, while $\tau(TT_n[V_5])=n$ and thus that the second inequality of Theorem~\ref{thm:bound-fvs} is tight.
Hence we conjecture the following.
\begin{conjecture}\label{conj:tight}
 For every positive integer $n$, there exists an oriented graph $D$ such that $\tau (D)=n$ and  $\inv(D) =2n$.
 \end{conjecture}


\medskip

A {\bf cycle packing} in a digraph is a set of vertex disjoint cycles.
The {\bf cycle packing number} of a digraph $D$, denoted by $\nu(D)$, is the maximum size of a cycle packing in $D$.
We have $\nu(D)\leq \tau(D)$ for every digraph $D$.
On the other hand, Reed et al.~\cite{RRST96} proved that there is a (minimum) function $f$ such that $\tau(D) \leq f(\nu(D))$ for every digraph $D$. With Theorem~\ref{thm:bound-fvs}, this implies $\inv(D) \leq  2\cdot f(\nu(D))$.

\begin{theorem}\label{thm:bound-nu}
There is a (minimum) function $g$ such that $\inv(D) \leq  g(\nu(D))$ for all oriented graph $D$ and $g\leq 2f$.
\end{theorem}

A natural question is then to determine this function $g$ or at least obtain good upper bounds on it.
Note that the upper bound on $f$ given by the proof of Reed et al.~\cite{RRST96} is huge (a multiply iterated exponential, where the number of iterations is also a multiply iterated exponential).
The only known value has been established by McCuaig~\cite{McCuaig91} who proved $f(1)=3$.
As noted in~\cite{RRST96}, the best lower bound on $f$ due to Alon (unpublished)  is $f(k)\geq k\log k$. It might be that $f(k)= O(k\log k)$. This would imply the following conjecture.

\begin{conjecture}\label{conj:klogk}
For all $k$, $g(k) =O(k\log k)$: 
there is an absolute constant $C$ such that $\inv(D) \leq C\cdot \nu(D) \log (\nu(D))$ for all oriented graph $D$.
\end{conjecture}

Note that for planar digraphs, combining results of Reed and Sheperd~\cite{ReSe96} and Goemans and Williamson~\cite{GoWi96}, we get $\tau(D) \leq 63\cdot \nu(D)$ for every planar digraph $D$.
This implies that $\tau(D) \leq 126\cdot \nu(D)$ for every planar digraph $D$ and so Conjecture~\ref{conj:klogk} holds for planar oriented graphs.

\medskip

Another natural question is whether or not the inequality $g\leq 2f$ is tight.
In Section~\ref{sec:intercyclic}, we show that it is not the case.
We show that $g(1)\leq 4$, while $f(1)=3$ as shown by McCuaig~\cite{McCuaig91}.
However we do not know if this bound $4$ on $g(1)$ is attained. Furthermore can we characterize the intercyclic digraphs with small inversion number~?

\begin{problem}
For any $k\in [4]$, can we characterize the intercyclic oriented graphs with inversion number $k$ ?
 \end{problem}



In contrast to Theorems~\ref{thm:bound-fvs} and \ref{thm:bound-nu}, the difference between $\inv$ and $\nu$, $\tau$, and $\tau'$ can be arbitrarily large as for every $k$, there are tournaments $T_k$ for which $\inv(T_k)=1$ and $\nu(T_k)= k$.
Consider for example the tournament $T_k$ obtained from three transitive tournaments $A$, $B$, $C$ of order $k$ by adding all arcs form $A$ to $B$, $B$ to $C$ and $C$ to $A$.
One easily sees that $\nu(T_k)=k$ and so $\tau'(T_k)\geq \tau(T_k)\geq k$; moreover $\Inv(T_k; A\cup B)$ is a transitive tournament, so $\inv(T_k)=1$.

\subsection{Maximum inversion number of an oriented graph of order $n$} 

For any positive integer $n$, let $\inv(n)=\max\{\inv(D) \mid D~\mbox{oriented graph of order}~n\}$. Since the inversion number is monotone (see Proposition~\ref{prop:monotone}), we have
$\inv(n)=\max\{\inv(T) \mid T~\mbox{tournament of order}~n\}$.

\begin{remark} 
 $\inv(n) \leq \inv(n-1)+1$ for all positive integer $n$.
\end{remark}
\begin{proof}
Let $T$ be a tournament of order $n$.
Pick a vertex $x$ of $T$. It is a sink in $D'=\Inv(T; N^+[x])$. 
So $\inv(D') =\inv(D'-x) \leq \inv(n-1)$ by Lemma~\ref{lem:reduc}. Hence $\inv(T) \leq \inv(n-1)+1$.
\end{proof}

Every oriented graph on at most two vertices is acyclic, so $\inv(1)=\inv(2)=0$.
Every tournament of order at most $4$ has a cycle arc-transversal of size at most $1$, so $\inv(3)=\inv(4)=1$.
As observed by Belkhechine et al.~\cite{BBBP}, every tournament of order at most $6$ has inversion number at most $2$. 

\begin{equation}\label{up-easy}
\inv(n) \leq n-4  ~~~~\mbox{for all}~n\geq 6.
\end{equation}

Moreover, Belkhechine et al.~\cite{BBBP10} observed that since there are $n!$ labelled transitive tournaments of order $n$, the number of labelled tournaments of order $n$ with inversion number less than $p$ is at most $n! 2^{n(p-1)}$, while there are $2^{\frac{n(n-1)}{2}}$ labelled tournaments of order $n$. So for all $n$ such that $2^{\frac{n(n-1)}{2}} >  n! 2^{n(p-1)}$, there is a tournament $T$ of order $n$ such that $\inv(T) \geq p$. Hence

\begin{equation}\label{down-easy}
\inv(n) \geq \frac{n-1}{2} - \log_2 n  ~~~~\mbox{for all}~n.
\end{equation}

However, it is believed that Equation~\eqref{down-easy} is not tight.

 \begin{conjecture}[Belkhechine et al.~\cite{BBBP}]
 $\inv(n) \geq \lfloor\frac{n-1}{2}\rfloor$.
 \end{conjecture}

 Furthermore, some explicit tournaments have been conjectured to have inversion number at least  $\lfloor\frac{n-1}{2}\rfloor$. 
 Let $Q_n$ be the tournament obtained from the transitive tournament by reversing the arcs of its unique directed hamiltonian path $(v_1,v_2, \ldots{},v_n)$.

 \begin{conjecture}[Belkhechine et al.~\cite{BBBP}]\label{conj:Q}
 $\inv(Q_n)= \lfloor\frac{n-1}{2}\rfloor$.
 \end{conjecture}

 A possible way to prove Conjecture~\ref{conj:Q} would be via augmentations.
 Let $D$ be an oriented graph and $z$ a vertex of $D$.
 The {\bf $z$-augmentation} of $D$ is the digraph, denoted by $\sigma(z,D)$,  obtained from $D$ by adding two new vertices $y$ and $x$, the arc $zy$, $yx$ and $xz$ and all the arcs from $\{x,y\}$ to $V(D)\setminus \{z\}$. 
We let $\sigma_i(z, D)$ be the $z$-augmentation of $D$ on which the vertices added are denoted by $x_i$ and $y_i$.
 
 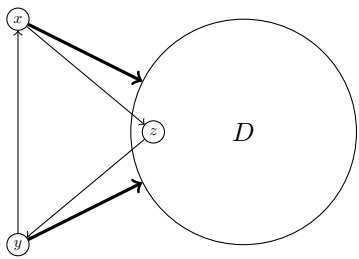
\begin{figure}[hbtp]
 \begin{center}
 \tikzset{invisible/.style={circle,minimum size=10pt,inner sep=0pt}}
\tikzset{blackvertex/.style={draw,circle,minimum size=5pt,inner sep=0pt, fill=black}}
\tikzset{edge/.style = {->,very thick}}
\tikzset{arc/.style = {->}}
\tikzset{box/.style={  rectangle, fill=white, rounded corners,  draw, minimum width=1.5cm, minimum height=1cm, inner sep=0pt}}
\tikzset{circ/.style={ circle,fill=white, minimum width=3cm, draw}}
\tikzstyle{vertex}=[circle,draw, minimum size=14pt, scale=0.6, inner sep=0.5pt]

  \begin{tikzpicture}
\node[circ] (D) at (2,0)  {$D$};

\node[vertex] (x) at (-1, 1.5) {$x$};
\node[vertex] (y) at (-1, -1.5) {$y$};
\node[vertex] (z) at (0.8, 0) {$z$};


\draw[edge] (x) edge (D);

\draw[edge] (y) edge (D);

\draw[arc] (z) edge (y);

\draw[arc] (y) edge (x);

\draw[arc] (x) edge (z);

\end{tikzpicture}
\caption{The $z$-augmentation $\sigma(z,D)$ of a digraph $D$.}
 \end{center}
\end{figure}

Observe that $Q_n$ is isomorphic to $\sigma(v_1,Q_{n-2})$.
Moreover for every oriented graph $D$ and vertex $z$ of $D$, 
$\inv(\sigma(z,D)) \leq \inv(D) + 1$, because $\Inv(\sigma(z,D), \{y,z\}) = (y \ra x) \ra D$.

In Section~\ref{sec:augment}, we prove that if $D$ is an oriented graph  with $\inv(D)=1$, then $\inv(\sigma(z,D)) = 2$ for every $z\in V(D)$ (Lemma~\ref{lem:auginv=1}).
 In particular, $\inv(Q_5)=2$.
 
\begin{figure}[hbtp]
 \begin{center}
 \tikzset{invisible/.style={circle,minimum size=10pt,inner sep=0pt}}
\tikzset{blackvertex/.style={draw,circle,minimum size=5pt,inner sep=0pt, fill=black}}
\tikzset{edge/.style = {->,very thick}}
\tikzset{arc/.style = {->}}
\tikzset{box/.style={  rectangle, fill=white, rounded corners,  draw, minimum width=1.5cm, minimum height=1cm, inner sep=0pt}}
\tikzset{circ/.style={ circle,fill=white, minimum width=3cm, draw}}
\tikzstyle{vertex}=[circle,draw, minimum size=14pt, scale=0.6, inner sep=0.5pt]

  \begin{tikzpicture}[scale=0.8]
 \node (a) at (0:1cm) [vertex] {$a$};
  \node (b) at (-120:1cm) [vertex] {$b$};
  \node (c) at (+120:1cm) [vertex] {$c$};

\node (x) at (2.5,1.2) [vertex] {$x_1$};
  \node (y) at (2.5,-1.2) [vertex] {$y_1$};
 
 \node (x') at (4,0.8) [vertex] {$x_2$};
  \node (y') at (4,-0.8) [vertex] {$y_2$};

\draw (0,0) circle (1.5);

 \node at (1.75,-1.5) {$H_1$};

 \draw [rounded corners] (-2,-2) rectangle (3,2);

  \draw [->, line width=0.02cm] (a) to (b);
\draw [->, line width=0.02cm] (b) to (c);
 \draw [->, line width=0.02cm] (c) to (a);
 
\draw [->, line width=0.02cm] (y) to (x);
\draw [->, line width=0.02cm] (y') to (x');

\draw [->, line width=0.02cm] (x) to (a);
\draw [->, line width=0.02cm] (a) to (y);
\draw [->, line width=0.02cm] (x') to (a);
\draw [->, line width=0.02cm] (a) to (y');

\draw [->, line width=0.04cm] (x) to (1.3,0.7);
\draw [->, line width=0.04cm] (y) to (1.3,-0.7);

\draw [->, line width=0.04cm] (x') to (3,0.8);
\draw [->, line width=0.04cm] (y') to (3,-0.8);

\end{tikzpicture}
\caption{The digraph $H_2$.}\label{fig:H'}
 \end{center}
\end{figure}
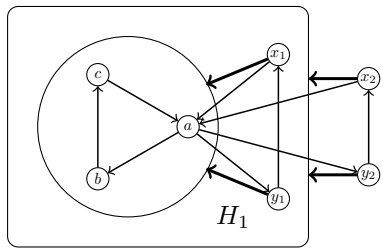

 Unfortunately, for larger values of $\inv(D)$, it is not true that $\inv(\sigma(z,D)) = \inv(D)+1$ for every $z\in V(D)$.
 For example, take the directed $3$-cycle $\vec{C_3}$ with vertex set $\{a,b,c\}$ and consider $H_1 = \sigma_1(a, \vec{C_3})$, 
and $H_2=  \sigma_2(a, H_1)$. See Figure~\ref{fig:H'}.
By  Lemma~\ref{lem:auginv=1}, we have $\inv(H_1)=2$ but $\inv(H_2)=2$ as
$(\{y_1, y_2, b\}, \{y_1, y_2, a, b\})$ is a decycling family of $H_2$. 
 
 However, we prove in Theorem~\ref{thm:2augment} that if $\inv(D)=1$, then  $\inv (\sigma_1(x_2,\sigma_2(z, D))) =3$ for every $z\in V(D)$. This directly implies $\inv(Q_7) = 3$.

 \subsection{Complexity of computing the inversion number}
 
We also consider the complexity of computing the inversion number of an oriented graph and the following associated problem.


\noindent {\sc $k$-Inversion}.\\
\underline{Input}: An oriented graph $D$. \\
\underline{Question}: $\inv(D)\leq k$ ? 

\medskip

We also study the complexity of the restriction of this problem to tournaments.


\noindent {\sc $k$-Tournament-Inversion}.\\
\underline{Input}: A tournament. \\
\underline{Question}: $\inv(T)\leq k$ ? 

 \smallskip

Note that {\sc $0$-Inversion} is equivalent to deciding whether an oriented graph $D$ is acyclic.
This can be done in $O(|V(D)|^2)$ time.
 
 Let $k$ be a positive integer.
 A tournament $T$ is {\bf $k$-inversion-critical} if $\inv(T) =k$ and $\inv(T-x) < k$ for all $x\in V(T)$.
We denote by ${\cal IC}_k$ the set of $k$-inversion-critical tournaments.
Observe that a tournament $T$ has inversion number at least $k$ if and only if $T$ has a subtournament in ${\cal IC}_k \cup {\cal IC}_{k+1}$ (by Lemma~\ref{lem:recur}).

\begin{theorem}[Belkhechine et al.~\cite{BBBP10}]\label{thm:crit-finite}
For  any positive integer $k$, the set ${\cal IC}_k$ is finite.
\end{theorem}

 Checking whether the given tournament $T$ contains $I$ for every element $I$ in ${\cal IC}_{k+1}\cup {\cal IC}_{k}$, one can decide  whether $\inv(T)\geq k$ in $O(|V(T)|^{\max
\{m_{k+1}, m_k\}})$ time, where $m_{k}$ is the maximum order of an element of ${\cal IC}_{k}$.
 
\begin{corollary}\label{cor:poly-tour}
For any non-negative integer $k$, {\sc $k$-Tournament-Inversion} is polynomial-time solvable. 
\end{corollary}  

The proof of Theorem~\ref{thm:crit-finite} neither explicitly describes  ${\cal IC}_{k}$ nor gives upper bound on $m_k$.
So the degree of the polynomial in Corollary~\ref{cor:poly-tour} is unknown.
This leaves open the following questions.
\begin{problem}
Explicitly describe ${\cal IC}_{k}$ or at least find an upper bound on $m_k$. \\
What is the minimum real number $r_k$ such that {\sc $k$-Tournament-Inversion} can be solved in $O(|V(T)|^{r_{k}})$ time ?
\end{problem}

\begin{figure}[hbtp]
\begin{center}
\tikzstyle{vertex}=[circle,draw, minimum size=12pt, scale=0.8, inner sep=0.5pt]

\begin{tikzpicture}[scale=0.8]
 \node (v1) at (90:1cm) [vertex] {$v_1$};
  \node (v2) at (-30:1cm) [vertex] {$v_2$};
  \node (v3) at (-150:1cm) [vertex] {$v_3$};

\begin{scope}[xshift=4.5 cm]
\node (v4) at (90:1cm) [vertex] {$v_4$};
  \node (v5) at (-30:1cm) [vertex] {$v_5$};
  \node (v6) at (-150:1cm) [vertex] {$v_6$};       
      \end{scope}

\draw (0,0) circle (1.5);
\draw (4.5,0) circle (1.5);

 \node at (2.25,-1.5) {$A_6$};
 
  \draw [->, line width=0.03cm] (v1) to (v2);
\draw [->, line width=0.03cm] (v2) to (v3);
 \draw [->, line width=0.03cm] (v3) to (v1);
 
\draw [->, line width=0.03cm] (v4) to (v5);
 \draw [->, line width=0.03cm] (v5) to (v6);
 \draw [->, line width=0.03cm] (v6) to (v4);

 \draw [->, line width=0.06cm] (1.5,0) to (3,0);

\end{tikzpicture}
\hfill
\begin{tikzpicture}[scale=0.8]
 \node (v3) at (90:1cm) [vertex] {$v_3$};
  \node (v4) at (-30:1cm) [vertex] {$v_4$};
  \node (v2) at (-150:1cm) [vertex] {$v_2$};

\node (v5) at (3,1) [vertex] {$v_5$};
  \node (v1) at (3,-1) [vertex] {$v_1$};

\draw (0,0) circle (1.5);

 \node at (1.75,-1.5) {$D_5$};
 
  \draw [->, line width=0.03cm] (v3) to (v4);
\draw [->, line width=0.03cm] (v4) to (v2);
 \draw [->, line width=0.03cm] (v2) to (v3);
 
\draw [->, line width=0.03cm] (v5) to (v1);

\draw [->, line width=0.03cm] (v1) to (1.4,-0.5);
\draw [<-, line width=0.03cm] (v5) to (1.4,0.5);

\end{tikzpicture}

\medskip

\begin{tikzpicture}[scale=0.8]
 \node (v1) at (0:1.8cm) [vertex] {$v_1$};
  \node (v2) at (-60:1.8cm) [vertex] {$v_2$};
  \node (v3) at (-120:1.8cm) [vertex] {$v_3$};
\node (v4) at (180:1.8cm) [vertex] {$v_4$};
  \node (v5) at (120:1.8cm) [vertex] {$v_5$};
  \node (v6) at (60:1.8cm) [vertex] {$v_6$};

 \node at (0,-2.3) {$B_6$};
 
  \draw [->, line width=0.03cm] (v1) to (v2);
\draw [->, line width=0.03cm] (v1) to (v3);
\draw [<-, line width=0.03cm] (v1) to (v4);
\draw [<-, line width=0.03cm] (v1) to (v5);
\draw [->, line width=0.03cm] (v1) to (v6);

\draw [->, line width=0.03cm] (v2) to (v3);
\draw [->, line width=0.03cm] (v2) to (v4);
\draw [<-, line width=0.03cm] (v2) to (v5);
\draw [<-, line width=0.03cm] (v2) to (v6);

 \draw [->, line width=0.03cm] (v3) to (v4);
  \draw [->, line width=0.03cm] (v3) to (v5);
  \draw [->, line width=0.03cm] (v3) to (v6);
  
\draw [<-, line width=0.03cm] (v4) to (v5);
 \draw [->, line width=0.03cm] (v4) to (v6);

 \draw [<-, line width=0.03cm] (v5) to (v6);

\end{tikzpicture}
\hfill
\begin{tikzpicture}[scale=0.8]

\node (v1) at  (18:1.8cm) [vertex] {$v_1$};
 \node (v2) at (-54:1.8cm) [vertex] {$v_4$};
 \node (v3) at (-126:1.8cm) [vertex] {$v_2$};
  \node (v4) at (162:1.8cm) [vertex] {$v_5$};
   \node (v5) at (90:1.8cm) [vertex] {$v_3$};

 \draw [->, line width=0.03cm] (v1) to (v2);
\draw [->, line width=0.03cm] (v1) to (v3);
 \draw [->, line width=0.03cm] (v2) to (v3);
\draw [->, line width=0.03cm] (v2) to (v4);
 \draw [->, line width=0.03cm] (v3) to (v4);
\draw [->, line width=0.03cm] (v3) to (v5);
 \draw [->, line width=0.03cm] (v4) to (v5);
\draw [->, line width=0.03cm] (v4) to (v1);
 \draw [->, line width=0.03cm] (v5) to (v1);
\draw [->, line width=0.03cm] (v5) to (v2);

 \node at (0,-2.2) {$R_5$};
\end{tikzpicture}
\hfill
\begin{tikzpicture}[scale=0.8]
 \node (v1) at (0.6,0) [vertex] {$v_1$};
  \node (v2) at (2.2,0) [vertex] {$v_2$};
  \node (v3) at (3.8,0) [vertex] {$v_3$};
 \node (v4) at (5.4,0) [vertex] {$v_4$};
  \node (v5) at (3,1) [vertex] {$v_5$};

 \node at (3,-2) {$V_5$};
 
  \draw [->, line width=0.03cm] (v1) to (v2);
\draw [->, line width=0.03cm] (v2) to (v3);
 \draw [->, line width=0.03cm] (v3) to (v4);
 \draw [->, line width=0.03cm] (v1) to [in=-150, out=-30] (v3);
 \draw [->, line width=0.03cm] (v2) to [in=-150, out=-30] (v4);
  \draw [->, line width=0.03cm] (v1) to [in=-135, out=-45] (v4);

\draw [<-, line width=0.03cm] (v1) to (v5);
 \draw [->, line width=0.03cm] (v2) to (v5);
 \draw [<-, line width=0.03cm] (v3) to (v5);
 \draw [->, line width=0.03cm] (v4) to (v5);
\end{tikzpicture}

\caption{The $2$-inversion-critical tournaments}\label{fig:2IC}
\end{center}
\end{figure}
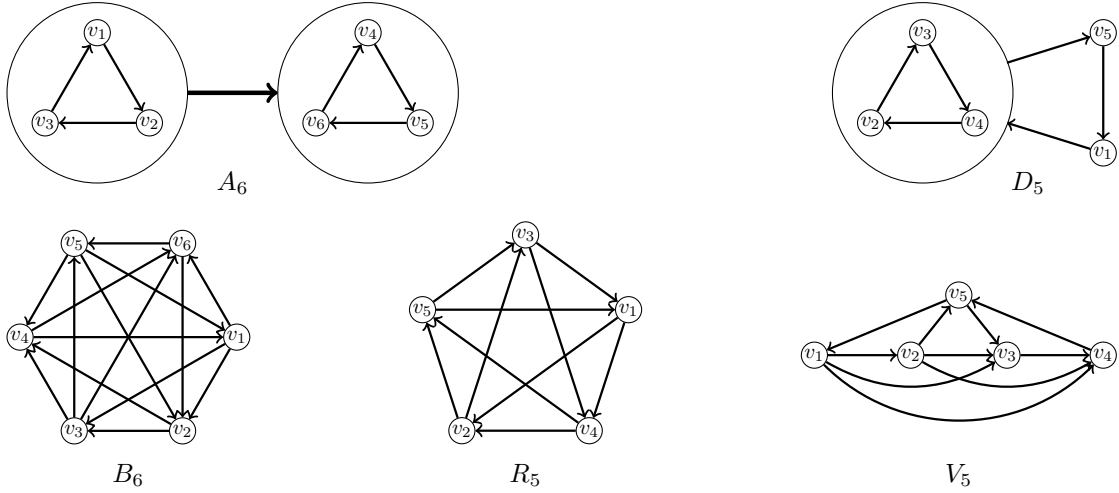

As observed in~\cite{BBBP10}, ${\cal IC}_{1}=\{\vec{C_3}\}$, so $m_1=3$.
This implies that {\sc $0$-Tournament-Inversion} can be done in $O(n^3)$. However, deciding whether a tournament is acyclic can be solved in $O(n^2)$-time.
Belkhechine et al.~\cite{BBBP10} also proved that ${\cal IC}_{2}=\{A_6, B_6, D_5, T_5, V_5\}$
where $A_6 = TT_2[\vec{C_3}] = \Inv(TT_6; (\{v_1, v_3\}, \{v_4, v_6\}))$,
$B_6= \Inv(TT_6; (\{v_1, v_4, v_5\}, \{v_2, v_5, v_6\}))$,
$D_5= \Inv(TT_5; (\{v_2,v_4\}, \{v_1,v_5\}))$,
$R_5= \Inv(TT_5; (\{v_1,v_3,v_5\}, \{v_2,v_4\}))$,
and $V_5= \Inv(TT_5; (\{v_1,v_5\}, \{v_3,v_5\}))$.
See Figure~\ref{fig:2IC}.

Hence $m_2=6$, so  {\sc $1$-Tournament-Inversion} can be solved in $O(n^6)$-time.
This is not optimal: we show in Subsection~\ref{subsec-poly} that it can be solved in $O(n^3)$-time, and that {\sc $2$-Tournament-Inversion} can be solved in $O(n^6)$-time.

\medskip
There is no upper bound on $m_k$ so far. 
Hence  since the inversion number of a tournament can be linear in its order (See e.g. tournament $T_k$ described at the end of the introduction), Theorem~\ref{thm:crit-finite} does not imply that one can compute the inversion number of a tournament in polynomial time.
In fact, we believe that it cannot be calculated in polynomial time.
 
\begin{conjecture}
Given a tournament and an integer $k$, deciding whether $\inv(T)=k$ is NP-complete. 
 \end{conjecture}

\medskip

In contrast to Corollary~\ref{cor:poly-tour}, we show in Subsection~\ref{subsec:NP} that {\sc $1$-Inversion} is NP-complete.
Note that together with Conjecture~\ref{conj:dijoin}, this would imply that {\sc $k$-Inversion} is NP-complete for every positive integer $k$.

  \begin{conjecture}
  {\sc $k$-Inversion} is NP-complete for all positive integer $k$.
    \end{conjecture}
    
As we proved  Conjecture~\ref{conj:dijoin}. when $\inv(L)=\inv(R)=1$, we get that {\sc $2$-Inversion} is NP-complete. 
    
\medskip    

Because of its relations with $\tau'$, $\tau$, and $\nu$, (see Subsection~\ref{subsec:versus}), it is natural to ask about the complexity of computing the inversion number when restricted to oriented graphs (tournaments) for which one of these parameters is bounded.
Recall that $\inv(D)=0$ if and only if $D$ is acyclic, so if and only if  $\tau'(D)= \tau(D) = \nu(D)=0$.

\begin{problem}
Let $k$ be a positive integer and $\gamma$ be a parameter in $\{\tau', \tau, \nu\}$.
What is the complexity of computing the inversion number of an oriented graph (tournament) $D$ with $\gamma(D) \leq k$ ?
\end{problem}

Conversely, it is also natural to ask about the complexity of  computing any of $\tau'$, $\tau$, and $\nu$, when restricted to oriented graphs  with bounded inversion number. In Subsection~\ref{subsec:bounded-inv}, we show that computing any of these parameters is NP-hard even for oriented graphs with inversion number $1$.
However, the question remains open when we restrict to tournaments.

\begin{problem}
Let $k$ be a positive integer and $\gamma$ be a parameter in $\{\tau', \tau, \nu\}$.
What is the complexity of computing $\gamma(T)$ for a tournament $T$ with $\inv(T) \leq k$ ?
\end{problem}

\section{Properties of the inversion number}\label{sec:easy}

In this section, we establish easy properties of the inversion number and deduce from them Theorem~\ref{thm:bound-fvs} and the fact that Conjecture~\ref{conj:dijoin} is equivalent to its restriction to tournaments.

The inversion number is monotone :

\begin{proposition}\label{prop:monotone}
If  $D'$ is a subdigraph of an oriented graph $D$, then $\inv(D')\leq \inv(D)$.
\end{proposition}
\begin{proof}
Let  $D'$ be a subdigraph of $D$.  
If $(X_i)_{i\in I}$ is  a decycling family of $D$, then 
$(X_i\cap V(D'))_{i\in I}$ is a decycling family of $D'$.
\end{proof}

 \begin{lemma}\label{lem:reduc}
 Let $D$ be an oriented graph.
 If $D$ has a source (a sink) $x$, then $\inv(D) = \inv(D-x)$.
 \end{lemma}
 \begin{proof}
 Every decycling family of $D-x$  is also a decycling family of $D$ since adding a source (sink) to an acyclic digraph results in an acyclic digraph.
 \end{proof}

\begin{lemma}\label{lem:recur}
Let $D$ be an oriented graph and let $x$ be a vertex of $D$.
Then $\inv(D) \leq \inv(D-x)+2$.
\end{lemma}
\begin{proof}
Let $N^+[x]$ be the closed out-neighbourhood of $x$, that is $\{x\}\cup N^+(x)$.
Observe that $D'=\Inv(D; (N^+[x], N^+(x)))$ is the oriented graph obtained from $D$ by reversing the arc between $x$ and its out-neighbours. Hence $x$ is a sink in $D'$ and $D'-x=D-x$. Thus, by Lemma~\ref{lem:reduc}, 
$\inv(D) \leq \inv(D')+2\leq  \inv(D-x)+2$.
\end{proof}

\begin{proof}[Proof of Theorem~\ref{thm:bound-fvs}]
As observed in the introduction, if $F$ is a minimum cycle arc-transversal, then the family of sets of end-vertices of arcs of $F$ is a decycling family.
So  $\inv(D) \leq  \tau' (D)$.

Let  $S=\{x_1, \dots , x_k\}$ be a cycle transversal with $k=\tau(D)$.
Lemma~\ref{lem:recur} and a direct induction imply 
$\inv(D) \leq \inv(D-\{x_1, \dots , x_i\}) +2i$  for all $i\in [k]$.
Hence   $\inv(D) \leq \inv(D-S) +2k$. But, since $S$ is a cycle transversal, $D-S$ is acyclic, so
$\inv(D-S) =0$. Hence $\inv(D)\leq 2k= 2\tau(D)$.
\end{proof}

Let $D$ be an oriented graph. An {\bf extension} of $D$ is any tournament $T$ such that $V(D)=V(T)$ and $A(D)\subseteq A(T)$.
\begin{lemma}\label{lem:extension}
Let $D$ be an oriented graph. There is an extension $T$ of $D$ such that $\inv(T)=\inv(D)$.
\end{lemma}
\begin{proof}
Set $p=\inv(D)$ and let $(X_i)_{i\in [p]}$ be a decycling family of $D$.
Then $D^*=\Inv(D; (X_i)_{i\in [p]})$ is acyclic and so admits an acyclic ordering $(v_1, \dots , v_n)$.

Let $T$ be the extension of $D$ constructed as follows: 
For every $1\leq k< \ell\leq n$ such that $v_kv_{\ell}\notin A(D^*)$, let $n(k,\ell)$ be the number of $X_i$, $i\in [p]$, such that $\{v_k,v_\ell\}\subseteq X_i$.
If $n(k,\ell)$ is even then the arc $v_kv_{\ell}$ is added to $A(T)$, and if $n(k,\ell)$ is odd then the arc $v_{\ell}v_k$ is added to $A(T)$.
Note that in the first case, $v_kv_{\ell}$ is reversed an even number of times by $(X_i)_{i\in [p]}$, and in the second $v_{\ell}v_k$ is reversed an odd number of times by $(X_i)_{i\in [p]}$. Thus, in both cases, $v_kv_{\ell}\in A(\Inv(T;  (X_i)_{i\in [p]}))$.
Consequently, $(v_1, \dots , v_n)$ is also an acyclic ordering of  $\Inv(T;  (X_i)_{i\in [p]})$.
Hence $\inv(T)\leq \inv(D)$, and so, by Proposition~\ref{prop:monotone}, $\inv(T)=\inv(D)$.
\end{proof}

\begin{proposition}\label{prop:equiv}
Conjecture~\ref{conj:dijoin} is equivalent to its restriction to tournaments.
\end{proposition}
\begin{proof}
Suppose there are oriented graphs $L,R$ that form a counterexample to Conjecture~\ref{conj:dijoin}, that is such that $\inv(L\ra R) < \inv(L) + \inv(R)$. By Lemma~\ref{lem:extension}, there is an extension $T$ of $L\ra R$ such that $\inv(T) = \inv(L\ra R)$ and let $T_L=T\langle V(L)\rangle$ and  $T_R=T\langle V(R)\rangle$. We have $T=T_L\ra T_R$ and by Proposition~\ref{prop:monotone}, $\inv(L)\leq \inv(T_L)$ and  $\inv(R)\leq \inv(T_R)$. Hence $\inv(T) < \inv(T_L) + \inv (T_R)$, so $T_L$ and $T_R$ are two tournaments that form a counterexample to Conjecture~\ref{conj:dijoin}.
\end{proof}

 \section{Inversion number of dijoins of oriented graphs}\label{sec:dijoin}

In this section, we give some evidence for Conjecture~\ref{conj:dijoin} to be true. We prove that it holds when $\inv(L)$ and $\inv(R)$ are small.


 \begin{proposition}\label{prop:1>1}
  Let \(L\) and \(R\) be two oriented graphs.
 If $\inv(L), \inv(R)\geq 1$, then $\inv(L\ra R)\geq 2$.
 \end{proposition} 
 \begin{proof}
 Assume $\inv(L), \inv(R)\geq 1$. Then $L$ and $R$ are not acyclic, so let $C_L$ and $C_R$ be directed cycles in $L$ and $R$ respectively.
 Assume for a contradiction that there is a set $X$ such that inverting $X$ in $L\ra R$ results in an acyclic digraph $D'$.
 There must be an arc $xy$ in $A(C_L)$ such that $x\in X$ and $y\notin X$, and there must be $z\in X\cap V(C_R)$.
 But then $(x, y, z, x)$ is a directed cycle in $D'$, a contradiction.
 \end{proof}

 Propositions~\ref{prop:D1>D2} and~\ref{prop:1>1} directly imply that Conjecture~\ref{conj:dijoin} holds when  $\inv(L) =  \inv(R) = 1$.
 \begin{corollary}\label{cor:1+1}
  Let \(L\) and \(R\) be two oriented graphs.
 If $\inv(L) = \inv(R) = 1$, then $\inv(L\ra R) = 2$.
 \end{corollary}

Further than Proposition~\ref{prop:1>1}, the following result gives some property of a minimum decycling family of $L\rightarrow R$ when  $\inv(L) = \inv(R) = 1$.

\begin{theorem}\label{thm:1>1}
 Let $D=(L\rightarrow R)$, where $L$ and $R$ are two oriented graphs with $\inv(L) = \inv(R) = 1$. Then, for any decycling family $(X_1, X_2)$ of $D$, either $X_1\subset V(L), X_2\subset V(R)$ or $X_1\subset V(R), X_2\subset V(L)$.
 \end{theorem}
  
 \begin{proof}
Let $(X_1, X_2)$ be a decycling family of $D$ and let $D^*$ be the acyclic digraph obtained after inverting $X_1$ and $X_2$ (in symbols $D^*=\Inv(D; (X_1,X_2))$).

Let us define some sets. See Figure \ref{fig:gstar}.

\begin{itemize}
\item For $i\in [2]$, $X_i^L = X_i\cap V(L)$ and $X_i^R = X_i\cap V(R)$. 
\item $Z^L = V(L)\setminus (X^L_1\cup X^L_2)$ and  $Z^R = V(R)\setminus (X^R_1\cup X^R_2)$.
\item $X^L_{12}=X^L_1\cap X^L_2$ and $X^R_{12}=X^R_1\cap X^R_2$.
\item for $\{i,j\}=\{1,2\}$, $X_{i-j}^L=(X^L_i\setminus X^L_j)$ and  $X_{i-j}^R=(X^R_i\setminus X^R_j)$.
 \end{itemize}
 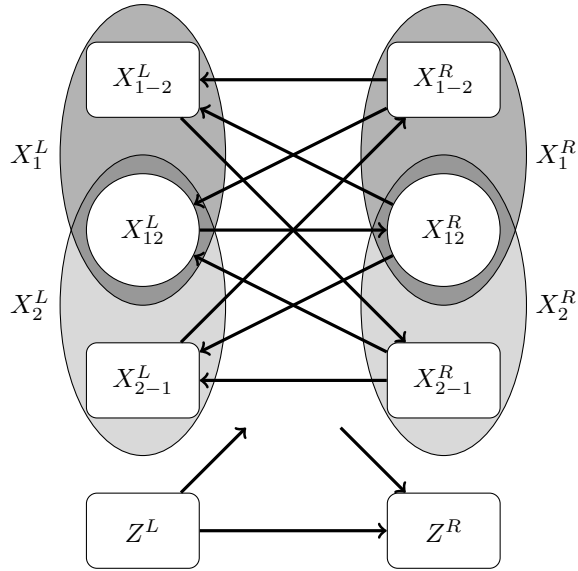
\begin{figure}[!ht]
  	\centering
\begin{tikzpicture}
\tikzset{invisible/.style={circle,minimum size=10pt,inner sep=0pt}}
\tikzset{blackvertex/.style={draw,circle,minimum size=5pt,inner sep=0pt, fill=black}}
\tikzset{edge/.style = {->,very thick}}
\tikzset{box/.style={  rectangle, fill=white, rounded corners,  draw, minimum width=1.5cm, minimum height=1cm, inner sep=0pt}}
\tikzset{circ/.style={ circle,fill=white, minimum width=1.5cm, draw}}

\begin{scope}[blend group=multiply]
    \draw[fill=gray!60!white] (-2,1) ellipse (1.1cm and 2cm);
    \draw[fill=gray!30!white] (-2,-1) ellipse (1.1cm and 2cm);
   	\draw[fill=gray!60!white] (2,1) ellipse (1.1cm and 2cm);
    \draw[fill=gray!30!white] (2,-1) ellipse (1.1cm and 2cm);  
\end{scope}

\node at (-3.5,1) {$X_1^L$};
\node at (-3.5,-1) {$X_2^L$};
\node at (3.5,1) {$X_1^R$};
\node at (3.5,-1) {$X_2^R$};

\node[box] (v1) at (-2,2) {$X_{1-2}^L$};

\node[circ] (v2) at (-2,0) {$X_{12}^L$};

\node[box] (v3) at (-2,-2) {$X_{2-1}^L$};

\node[box] (v4) at (2,2) {$X_{1-2}^R$};

\node[circ] (v5) at (2,0)  {$X_{12}^R$};

\node[box] (v6) at (2,-2) {$X_{2-1}^R$};
 
\node[box] (zd) at (-2,-4) {$Z^L$};
\node[invisible](n1) at (-0.5, -2.5) {};
\draw[edge] (zd) edge (n1){};

\node[box] (zf) at (2,-4) {$Z^R$};
\node[invisible](n1) at (0.5, -2.5) {};
\draw[edge] (n1) edge (zf){};



\draw[edge] (v1) edge (v6);

\draw[edge] (v2) edge (v5);

\draw[edge] (v3) edge (v4);

\draw[edge] (v4) edge (v1);
\draw[edge] (v4) edge (v2);

\draw[edge] (v5) edge (v1);
\draw[edge] (v5) edge (v3);

\draw[edge] (v6) edge (v2);
\draw[edge] (v6) edge (v3);

\draw[edge] (zd) edge (zf); 

\end{tikzpicture}

  	\caption{The oriented graph $D^*$} \label{fig:gstar}
  \end{figure}  
 
 Observe that at least one of the sets $X^L_{1-2}, X^R_{2-1}, X^L_{2-1}$ and $X^R_{1-2}$ must be empty, otherwise $D^*$ is not acyclic. 
 By symmetry, we may assume that it is $X^R_{1-2}$ or  $X^R_{2-1}$.
 Observe moreover that $X^R_{1-2} \cup X^R_{2-1} \neq \emptyset$ for otherwise $X_1^R=X_2^R=X^R_{12}$ and $D^*\langle V(R)\rangle =R$ is not acyclic.

  \medskip
  
  Assume first that $X^R_{1-2} =\emptyset$ and so $X^R_{2-1}\neq\emptyset$. 
  
  Suppose for a contradiction that $X^R_{12}\neq \emptyset$ and let $a\in X^R_{2-1}, b\in X^R_{12}$. Let $C$ be a directed cycle in $L$. Note that $V(C)$ cannot be contained in one of the sets $X^L_{1-2}, X^L_{12}$, $X^L_{2-1}$ or $Z^L$. If $V(C)\cap Z^L\neq \emptyset$, there is an arc $cd\in A(L)$ such that $c\in X^L_{1-2}\cup X^L_{12}\cup X^L_{2-1}$ and $d\in Z^L$. Then, either $(c, d, a, c)$ or $(c, d, b, c)$ is a directed cycle in $D^*$, a contradiction.  Thus, $V(C)\subseteq X^L_{1-2}\cup X^L_{12}\cup X^L_{2-1}$. If $V(C)\cap X^L_{12}\neq \emptyset$, then there is an arc $cd\in A(L)$ such that $c\in X^L_{12}$ and $d\in X^L_{1-2}\cup X^L_{2-1}$ which means that $dc\in A(D^*)$ and $(d, c, b, d)$ is a directed cycle in $D^*$, a contradiction. Hence $V(C)\subseteq X^L_{1-2}\cup X^L_{2-1}$ and there exists an arc $cd\in A(L)$ such that $c\in X^L_{2-1}, d\in X^L_{1-2}$ and $(c, d, a, c)$ is a directed cycle in $D^*$, a contradiction.
  
   Therefore $X^R_{12}=\emptyset$ and every directed cycle of $R$ has its vertices in $X^R_{2-1}\cup Z^R$. Then, there is an arc $ea\in A(R)$ with $a\in X^R_{2-1}$ and $e\in Z^R$. Note that, in this case, $ea\in A(D^*)$ and $(e, a, c, e)$ is a directed cycle in $D^*$ for any $c\in X^L_{12}\cup X^L_{2-1}$. Thus, $X^L_{12}=X^L_{2-1}=\emptyset$ and $X_1\subset V(L), X_2\subset V(R)$.

  \medskip
  If $X^R_{2-1}=\emptyset$, we can symmetrically apply the same arguments to conclude that $X_1\subset V(R)$ and $X_2\subset V(L)$. 
 \end{proof}

\begin{theorem}\label{thm:1+2=3}
 Let $L$ and $R$ be two oriented graphs. If $\inv(L) = 1$ and $\inv (R)=2$, then $\inv(L\rightarrow R) = 3$.
 \end{theorem}

 \begin{proof}
Let $D=(L\rightarrow R)$. By Propositions~\ref{prop:D1>D2} and \ref{prop:1>1}, we know that $2\leq \inv(D)\leq 3$. 

\medskip

Assume for a contradiction that $\inv(D) = 2$. 
Let $(X_1, X_2)$ be a decycling family of $D$ and let $D^*=\Inv(D; (X_1,X_2))$. 
Let $L^*=D^*\langle V(L)\rangle$ and $R^*=D^*\langle V(F)\rangle$. 
We define the sets $X^L_1$, $X^L_2$, $X^R_1$, $X^R_2$, $Z^L$, $Z^R$, $X^L_{12}$, $X^R_{12}$, 
$X^L_{1-2}$, $X^L_{2-1}$, $X^R_{1-2}$, and $X^R_{2-1}$ as in Theorem~\ref{thm:1>1}.
See Figure \ref{fig:gstar}. Note that each of these sets induces an acyclic digraph in $D^*$ and thus also in $D$.
For $i\in [2]$, let $D_i=\Inv(D;X_i)$, let $L_i=\Inv(L, X^L_i) = \Inv(L^*; X^L_{j-i})$ where $\{j\}=[2]\setminus \{i
\}$, and $R_i=\Inv(R, X^R_i) = \Inv(R^*; X^R_{j-i})$ where $\{j\}=[2]\setminus \{i\}$. 
Since $\inv(D)=2$, $\inv(D_1) = \inv(D_2) =1$.
Since $\inv(R)=2$, $R_1$ and $R_2$ are both non-acyclic, so $\inv(R_1) = \inv(R_2) =1$.
 
  \medskip
 
    \noindent{}{\bf Claim 1:}
 {\it    $X^L_i,X^R_i\neq\emptyset$ for all $i\in [2]$.}
    \begin{subproof}
   Since $\inv(R) = 2$, necessarily,  $X^R_1, X^R_2\neq\emptyset$.
   
 Suppose now that $X^L_i=\emptyset$ for some $i\in [2]$. Then $D_i = L \ra  R_i$. As $\inv(L)\geq 1$ and $\inv(R_i)\geq 1$, by Proposition~\ref{prop:1>1} 
 $\inv(D_i) \geq 2$, a contradiction.  \end{subproof}

    \noindent{}{\bf Claim 2:}
 {\it    $X^L_1\neq X^L_2$ and $X^R_1\neq X^R_2$.}
 \begin{subproof}
If $X_1^L=X_2^L$, then $L^*=L$, so $L^*$ is not acyclic, a contradiction.
Similarly, If $X_1^R=X_2^R$, then $R^*=R$, so $R^*$ is not acyclic, a contradiction.
 \end{subproof}
 
 In particular, Claim 2 implies that $X^L_{1-2} \cup X^L_{2-1} \neq \emptyset$.

  \medskip
  
  In the following, we denote by $A\leadsto B$ the fact that there is no arc from $B$ to $A$.
  
  \medskip
  
  Assume first that $X^R_{1-2} =\emptyset$. By Claim~1, $X^R_1\neq \emptyset$, so $X^R_{12}\neq \emptyset$ and by Claim~2 , $X^R_1\neq X^R_2$, so $X^R_{2-1}\neq\emptyset$.
  
   If $X^L_{2-1}\neq\emptyset$, then, in $D^*$,
$X^R_{2-1}\cup X^R_{12}\leadsto Z^R$ because $X^R_{2-1}\cup X^R_{12}\ra X^L_{2-1} \ra Z^R$. But then $R_1=\Inv(R^*; X^R_{2})$ would be acyclic, a contradiction. Thus, $X^L_{2-1}=\emptyset$.
 
 Then by Claims~1~and~2, we get  $X^L_{12}, X^L_{1-2}\neq \emptyset$. 
 Hence, as  $X^R_{12}\ra X^L_{1-2} \ra X^R_{2-1} \ra X^L_{12}\ra X^R_{12}$ in $D^*$, there is a directed cycle in $D^*$, a contradiction. Therefore $X^R_{1-2}\neq \emptyset$.
 
 \medskip
 In the same way, one shows that  $X^R_{2-1}\neq \emptyset$.
 As $X^R_{1-2} \ra X^L_{1-2} \ra X^R_{2-1}  \ra X^L_{2-1} \ra X^R_{1-2}$ in $D^*$, and $D^*$ is acyclic, one of $X^L_{1-2}$ and $X^L_{2-1}$ must be empty.
 Without loss of generality, we may assume $X^L_{1-2}=\emptyset$.
 
 Then   by Claims~1~and~2, we have $X^L_{12}, X^L_{2-1}\neq \emptyset$. 
 Furthermore $X^R_{12}=\emptyset$ because $X^R_{12} \ra X^L_{2-1} \ra X^R_{1-2}\ra X^L_{12} \ra X^R_{12}$ in $D^*$.
 Now  in $D^*$,
$X^R_{2-1}\leadsto X^R_{1-2} \cup Z^R$ because $X^R_{2-1}\ra X^L_{2-1} \ra X^R_{1-2} \cup Z^R$, and
$X^R_{1-2} \leadsto Z^R$ because $X^R_{1-2} \ra X^L_{12} \ra Z^R$.
Thus, in $D$, we also have $X^R_{2-1}\leadsto X^R_{1-2} \cup Z^R$ and $X^R_{1-2} \leadsto Z^R$.
So $R$ is acyclic, a contradiction to $\inv(R)\geq 2$. 
\medskip
  
    Therefore $\inv(D)\geq 3$.  So  $\inv(D)= 3$.  
 \end{proof}

 \begin{corollary}\label{cor:1+1=2}
Let $D$ be an oriented graph. Then $\inv(D) = 1$ if and only if $\inv(D\ra D)= 2$.
 \end{corollary}
 \begin{proof}
 Assume first that  $\inv(D) = 1$. Then by Corollary~\ref{cor:1+1}, $\inv(D\ra D)= 2$.\\
\indent Assume now that $\inv(D) \neq 1$. \\
 If $\inv(D) =0$, then $D$ is acyclic, and so is $D\ra D$. Hence $\inv(D\ra D)=0$.\\
 If $\inv(D) \geq 3$, then $\inv(D\ra D) \geq \inv (D)$ (by Proposition~\ref{prop:monotone} because $D$ is a subdigraph of $D\ra D$) and so $\inv(D\ra D)\geq 3$.\\
 If $\inv(D) = 2$, then $D$ contains a directed cycle $C$. Now $C\ra D$ is a subdigraph of  $D\ra D$, so by Proposition~\ref{prop:monotone}  $\inv(D\ra D) \geq \inv (C\ra D)$.
 Clearly, $\inv(C)=1$, thus, by Theorem~\ref{thm:1+2=3}, $\inv(C\ra D) = 3$ and so $\inv(D\ra D)\geq 3$.
  \end{proof}
 
 \subsection{Dijoin of oriented graphs with inversion number $2$}
  
  \begin{theorem}
  \label{thm:inv2and2}
  Let $L$ and $R$ be strong oriented graphs such that $\inv{}(L),\inv{}(R)\geq 2$. Then $\inv{}(L\rightarrow{}R)\geq 4$.
\end{theorem}
  \begin{proof}
  Assume for a contradiction that there are two strong oriented graphs $L$ and $R$ such that  $\inv{}(L),\inv{}(R)\geq 2$ and $\inv{}(L\rightarrow{}R)\leq 3$. By Lemma~\ref{lem:extension} and Proposition~\ref{prop:monotone}, we can assume that $L$ and $R$ are strong tournaments.

Hence $L$ contains $\vec{C_3}$. By Theorem~\ref{thm:1+2=3},  $\inv{}(\vec{C_3}\rightarrow{}R)\geq 3$. But $\vec{C_3}\rightarrow{}R$ is a subtournament of  $L\rightarrow{}R$. Thus, by Proposition~\ref{prop:monotone}, $\inv{}(L\rightarrow{}R)\geq 3$ and so $\inv{}(L\rightarrow{}R)=3$. Let $(X_1,X_2,X_3)$ be a decycling sequence of $D=L\rightarrow{}R$ and denote the resulting acyclic (transitive) tournament by $T$. We will use the following notation. Below and in the whole proof, whenever we use subscripts $i,j,k$ together we have $\{i,j,k\}=\{1,2,3\}$.

  \begin{itemize}
    
  \item $X^L_i=X_i\cap V(L), X^R_i= X_i\cap V(R)$ for all $i\in [3]$.
\item $Z^L=V(L)\setminus (X^L_1\cup X^L_2\cup X^L_3)$ and $Z^R=V(R)\setminus (X^R_1\cup X^R_2\cup X^R_3)$.
  \item $X^L_{123}=X^L_1\cap X^L_2\cap X^L_3$, $X^R_{123}=X^R_1\cap X^R_2\cap X^R_3$.
  \item  $X^L_{ij-k}=(X^L_i\cap X^L_j)\setminus X^L_k$ and $X^R_{ij-k}=(X^R_i\cap X^R_j)\setminus X^R_k$.
  \item $X^L_{i-jk}=X^L_i\setminus (X^L_j\cup X^L_k)$ and $X^R_{i-jk}=X^R_i\setminus (X^R_j\cup X^R_k)$.
    \end{itemize}

    For any two (possibly empty) sets $Q,W$, we write $Q\ra W$ to indicate that every $q\in Q$ has an arc to every $w\in W$. Unless otherwise specified, we are always referring to the arcs of $T$ below. When we refer to arcs of the original digraph we will use the notation $u\rad v$, whereas we use $u\ra v$ for arcs in $T$.\\

    \noindent{}{\bf Claim A:}
 {\it    $X^L_i,X^R_i\neq\emptyset$ for all $i\in [3]$.}
    \begin{subproof}
 Suppose w.l.o.g. that $X^R_1=\emptyset$ and let $D'=\Inv(D;X_1)$. Then $D'$ contains $\vec{C_3}\ra R$ as a subtournament since reversing $X^L_1$ does not make $L$ acyclic so there is still a directed 3-cycle (by Moon's theorem).
 \end{subproof}   

\noindent{}{\bf Claim B:}
{\it In $T$ the following holds, implying that at least one of the involved sets is empty (as $T$ is acyclic).
\begin{enumerate}
\item[\rm (a)] $X^R_{123}\ra X^L_{123}\ra X^R_{ij-k}\ra X^L_{ik-j}\ra X^R_{123}$.
\item[\rm (b)] $X^L_{ij-k}\ra X^R_{ij-k}\ra X^L_{ik-j}\ra X^R_{ik-j}\ra X^L_{ij-k}$.
      \end{enumerate}}
 \begin{subproof}
      This  follows from the fact that and arc of $D$ is inverted if and only if it belongs to an odd number of the sets $X_1,X_2,X_3$. 
 \end{subproof}

      \noindent{}{\bf Claim C:} {\it For all $i\neq j$, we have $X^L_i\neq X^L_j$ and $X^R_i\neq X^R_j$.}
 \begin{subproof}
      Suppose this is not true, then without loss of generality $X^L_3=X^L_2$ but this contradicts that $(X^L_1,X^L_2,X^L_3)$ is  a  decycling sequence of $L$ as inverting $X^L_2$ and $X^L_3$ leaves every arc unchanged and we have $\inv(L)\geq 2$.
\end{subproof}      

      Now we are ready to obtain a contradiction to the assumption that $(X_1,X_2,X_3)$ is a decycling sequence for $D=L\ra R$. We divide the proof into five cases. In order to increase readability, we will emphasize partial  conclusions in \fact{blue}, assumptions in \asm{orange}, and indicate consequences of assumptions in \sfact{red}.  \\

      \noindent{\bf Case 1:} $X^L_{i-jk}=\emptyset=X^R_{i-jk}$ for all $i,j,k$.\\
      
      By Claim C, at least two of the sets $X^L_{12-3},X^L_{13-2},X^L_{23-1}$ are non-empty and at least two of the sets $X^R_{12-3},X^R_{13-2},X^R_{23-1}$ are non-empty. Without loss of generality, \fact{$X^L_{12-3},X^L_{13-2}\neq\emptyset$}. Now Claim B~(b) implies that one of
      $X^R_{12-3},X^R_{13-2}$ must be empty. By interchanging the names of $X_2,X_3$ if necessary, we may assume that \fact{$X^R_{13-2}=\emptyset$} and hence, by Claim C, \fact{$X^R_{12-3},X^R_{23-1}\neq\emptyset$}. 
 By Claim B~(a), this implies $X^L_{23-1}=\emptyset$. 
      Now $X^R_{23-1}\ra X^L_{12-3}\ra X^R_{12-3}$, so 
      \fact{$X^R_{23-1}\ra X^R_{13-2}$}. As $X^L_{12-3}\ra X^R_{12-3}\ra X^L_{13-2}$, we must have \fact{$X^L_{12-3}\ra X^L_{13-2}$}.
      By Claim B~(a), $X^L_{123}\ra X^R_{12-3}\ra X^L_{13-2}\ra X^R_{123}\ra X^L_{123}$, so one of $X^L_{123}$ and $X^R_{123}$ is empty. W.l.o.g. we may assume \fact{$X^R_{123}=\emptyset$}. 
 As $R$ is strong and  $X^R_{23-1}$ dominates $X^R_{12-3}$ in $R$ (these arcs are reversed by $X_2$), we must have \fact{$Z^R\neq\emptyset$}.  Moreover the arcs incident to $Z^R$ are not reversed, so the set $Z^R$ has an out-neighbour in $X^R_{12-3}\cup X^R_{23-1}$. But $X^R_{12-3}\cup X^R_{23-1}\ra X^L_{13-2}\ra Z^R$ so $T$ has a directed 3-cycle, contradiction. This completes the proof of Case 1.\\

      \noindent{\bf Case 2:} Exactly one of $X^L_{1-23},X^L_{2-13},X^L_{3-12},X^R_{1-23},X^R_{2-13},X^R_{3-12}$ is non-empty.\\
      
      By reversing all arcs and switching the names of $L$ and $R$ if necessary, we may assume w.l.o.g that \fact{$X^L_{1-23}\neq\emptyset$}. As $X^R_2\neq X^R_3$ we have $X^R_{12-3}\cup X^R_{13-2}\neq\emptyset$. By symmetry, we can assume that \fact{$X^R_{12-3}\neq\emptyset$}.\\

      Suppose for a contradiction that
      \asm{$X^R_{23-1}=\emptyset$}. Then Claims A and  C imply \sfact{$X^R_{13-2}\neq\emptyset$}. Now, by Claim B~(b), one of $X^L_{12-3},X^L_{13-2}$ is empty. By symmetry, we can assume \sfact{$X^L_{13-2}=\emptyset$}. Now,  by Claim C, $X^L_2\neq X^L_3$, so  \sfact{$X^L_{12-3}\neq\emptyset$}. Note that
      $X^L_{12-3}\ra X^R_{12-3}\ra X^L_{1-23}$, thus
      \sfact{$X^L_{12-3}\ra X^L_{1-23}$} because $T$ is acyclic.  We also have \sfact{$X^L_{123}\ra X^L_{12-3}$} as $X^L_{123}\ra X^R_{13-2}\ra X^L_{12-3}$, and \sfact{$X^L_{12-3}\ra X^L_{23-1}$} as $X^L_{12-3}\ra X^R_{12-3}\ra X^L_{23-1}$. 
      This implies that in $L$ all arcs between $X^L_{12-3}$ and $X^L_{23-1}\cup X^L_{123}\cup X^L_{1-23}$ are entering $X^L_{12-3}$ (the arcs between $X^L_{123}$ and $X^L_{12-3}$ were reversed twice and those between $X^L_{1-23}\cup X^L_{23-1}$ and $X^L_{12-3}$ were reversed once). Hence, as $L$ is strong, we must have an arc $uz$ from $X^L_{12-3}$ to $Z^L$. But $Z^L\ra X^R_{13-2}\ra X^L_{12-3}$ so together with $uz$ we have a directed $3$-cycle in $T$, contradiction. Hence \fact{$X^R_{23-1}\neq\emptyset$}. \\

      Observe that \fact{$X^R_{12-3}\cup X^R_{13-2}\ra X^R_{23-1}$} as $X^R_{12-3}\cup X^R_{13-2}\ra X^L_{1-23}\ra X^R_{23-1}$.\\

If  \asm{$X^L_{12-3}\neq\emptyset$}, then $X^R_{23-1}\ra X^L_{12-3} \ra X^R_{12-3} \ra X^R_{23-1}$, a contradiction. So \fact{$X^L_{12-3}=\emptyset$}.
But $X_2^L\neq X_3^L$ by Claim C. Thus \fact{$X^L_{13-2}\neq\emptyset$}.
As $X^R_{23-1}\ra X^L_{13-2}\ra X^R_{123}$, we have \fact{ $X^R_{23-1}\ra X^R_{123}$}.
This implies that in $R$ all the arcs between $X^R_{23-1}$ and $X^R_{13-2}\cup X^R_{123}\cup X^R_{12-3}$ are leaving $X^R_{23-1}$. So as $R$ is strong there must be an arc in $R$ from $Z^R$ to $X^R_{23-1}$. This arc is not reversed, so it is also an arc in $T$. But since $X^R_{23-1}\ra X^L_{13-2}\ra Z^R$, this arc is in a directed $3$-cycle, a contradiction.
This completes Case 2.\\

      \noindent{}{\bf Case 3:} Exactly one of $X^L_{1-23},X^L_{2-13},X^L_{3-12}$ is non-empty and exactly one of  $X^R_{1-23},X^R_{2-13},X^R_{3-12}$ is non-empty.\\

      By symmetry we can assume \fact{$X^L_{1-23}\neq\emptyset$}. \\

      \noindent{}{\bf Subcase 3.1:} $X^R_{1-23}\neq\emptyset$.\\
  
By Claim~C, $X_2^L\neq X_3^L$, so one of $X^L_{12-3}$ and $X^L_{13-2}$ is non-empty.   By symmetry we may assume \fact{$X^L_{12-3}\neq\emptyset$}.\\
     
      Suppose \asm{$X^R_{12-3}\neq\emptyset$}. Then \sfact{$X^R_{23-1}=\emptyset$} as $X^L_{1-23}\ra X^R_{23-1}\ra X^L_{12-3}\ra X^R_{12-3}\ra X^L_{1-23}$, and \sfact{$X^L_{23-1}=\emptyset$} as
      $X^R_{1-23}\ra X^L_{12-3}\ra X^R_{12-3}\ra X^L_{23-1}\ra X^R_{1-23}$.
    
By Claim B~(b), one of  $X^L_{13-2},X^R_{13-2}$ is empty. By symmetry, we may assume \sfact{$X^R_{13-2}=\emptyset$}.

Observe that $V(R)\setminus Z^R= X^R_{123} \cup X^R_{12-3}\cup X^R_{1-23}$, so $V(R)\setminus Z^R \ra X^L_{1-23} \ra Z^R$, so $V(R)\setminus Z^R \ra Z^R$. But all the arcs incident to $Z^R$ are not inversed, so in $R$, there is no arc from $Z^R$ to  $V(R)\setminus Z^R$. Since $R$ is strong, \sfact{$Z^R=\emptyset$}. 

Now\sfact{ $X^R_{1-23} \ra X^R_{12-3}\cup X^R_{123}$} because $X^R_{1-23}\ra X^L_{12-3} \ra X^R_{12-3}\cup X^R_{123}$.
But all the arcs between $X^R_{1-23}$ and $X^R_{12-3}\cup X^R_{123}=V(R)\setminus X^R_{1-23}$ are inversed from $R$ to $T$.
Hence in $R$, no arcs leaves $X^R_{1-23}$ in $R$, a contradiction to $R$ being strong.

Hence \fact{$X^R_{12-3}=\emptyset$}. As $X^R_2\neq X^R_3$ this implies \fact{$X^R_{13-2}\neq\emptyset$}.\\

      Suppose that \asm{$X^R_{23-1}=\emptyset$}, then  \sfact{$X^R_{123}\neq\emptyset$} because $X^R_2\neq \emptyset$ by Claim A. Furthermore \sfact{$X^R_{13-2}\ra X^R_{123}$} as $X^R_{13-2}\ra X^L_{12-3}\ra X^R_{123}$, and  \sfact{$X^L_{12-3}\ra X^L_{1-23}$} as $X^L_{12-3}\ra X^R_{123}\ra X^L_{1-23}$. 
      This implies that \sfact{$X^L_{123}=\emptyset$} as $X^L_{123}\ra X^R_{13-2}\ra X^R_{123}\ra X^L_{123}$. 
      
       Since $L$ is strong, there must be an arc $uv$ leaving $X^L_{12-3}$ in $L$. But $v$ cannot be in  $X^L_{1-23}$ since all vertices of this set dominate $X^L_{12-3}$ in $L$. Moreover $v$ cannot be in $Z^L$ for otherwise $(u,v,w,u)$ would be a directed $3$-cycle in $T$  for any $w\in X^R_{1-23}$ since $Z^L\ra X^R_{1-23} \ra X^L_{12-3}$. 
Hence $v\in X^L_{13-2}\cup X^L_{23-1}$, so  \sfact{$X^L_{13-2}\cup X^L_{23-1}\neq\emptyset$}. 
      
      As $X^L_{13-2}\ra X^R_{13-2}\ra X^L_{23-1}\ra X^R_{1-23}\ra X^L_{13-2}$, precisely one of $X^L_{13-2},X^L_{23-1}$ is non-empty.

      If  $X^L_{13-2}\neq\emptyset$ and $X^L_{23-1}=\emptyset$, then $X^L_{13-2}\ra X^R_{13-2}\ra X^L_{1-23}\cup X^L_{12-3}$ implies that $X^L_{13-2}\ra X^L_{1-23}\cup X^L_{12-3}$. As $d^+_L(X^L_{13-2})>0$ there exists $z\in Z^L$ such that there is an arc $uz$ from $X^L_{13-2}$ to $Z^L$, but then $z\ra X^R_{1-23}\ra u\ra z$ is a contradiction. Hence \sfact{$X^L_{13-2}=\emptyset$} and \sfact{$X^L_{23-1}\neq\emptyset$}. Then \sfact{$X^L_{23-1}\ra X^L_1$} as $X^L_{23-1}\ra X^R_{1-23}\ra X^L_1$.

      Note that \sfact{$Z^L=\emptyset$} as every vertex in $V(L)\setminus Z^L$ has an in-neighbour in $V(R)$ in $T$, implying that there can be no arc from $V(L)\setminus Z^L$ to $Z^L$ in $L$. Thus $V(L)=X^L_{1-23}\cup X^L_{12-3}\cup X^L_{23-1}$ where each of these sets induces an acyclic subtournament of  $L$ and we have $X^L_{1-23}\rad X^L_{12-3}\rad X^L_{23-1}\rad X^L_{1-23}$ in $L$. But now inverting the set $X^L_{1-23}\cup X^L_{23-1}$ makes $L$ acyclic, a contradiction to $\inv(L)\geq 2$. Thus \fact{$X^R_{23-1}\neq \emptyset$.}\\
      
      Suppose  \asm{$X^L_{23-1}=\emptyset$}. As above  \sfact{$Z^L=\emptyset$}, so $V(L)=X^L_1$.
       As $X^L_{1-23}\ra X^R_{23-1}\ra X^L_{12-3}\cup X^L_{13-2}$ we have $X^L_{1-23}\ra X^L_{12-3}\cup X^L_{13-2}$.
       Thus, using $d^+_L(X^L_{1-23})>0$, we get \sfact{$X^L_{123}\neq\emptyset$}. As $X^R_{123}\ra X^L_{123}\ra X^R_{23-1}\ra X^L_{12-3}\ra X^R_{123}$, we have \sfact{$X^R_{123}=\emptyset$}. 
      Moreover \sfact{$X^R_1\ra X^R_{23-1}$} because $X^R_1\ra X^L_{1-23}\ra X^R_{23-1}$. We also have \sfact{$X^R_{1-23}\ra X^R_{13-2}$} as $X^R_{1-23}\ra X^L_{123}\ra X^R_{13-2}$.  
  Now $V(R)\setminus Z^R =X^R_{1-23} \cup X^R_{13-2} \cup X^R_{23-1} \ra  X^L_{12-3} \ra Z^R$. Thus $Z^R=\emptyset$ and  $V(R)=X^R_{1-23} \cup X^R_{13-2} \cup X^R_{23-1}$ where each of these sets induces an acyclic subtournament in $R$ and $X^R_{1-23}\rad X^R_{23-1}\rad X^R_{13-2}\rad X^R_{1-23}$ in $D$. But then inverting $X^R_{1-23}\cup X^R_{23-1}$ we make $R$ acyclic, a contradiction to $\inv(R)\geq 2$.
Thus \fact{$X^L_{23-1}\neq \emptyset$}.\\


Therefore  \fact{$X^L_{13-2}=\emptyset$} as
$X^L_{13-2}\ra X^R_{13-2}\ra X^L_{23-1}\ra X^R_{23-1}\ra X^L_{13-2}$. 
As $X^L_{1-23}\ra X^R_{23-1}\ra X^L_{12-3}$ we have \sfact{$X^L_{1-23}\ra X^L_{12-3}$}; as $X^L_{23-1}\ra X^R_{1-23}\ra X^L_{1-23}\cup X^L_{12-3}$ we have \sfact{$X^L_{23-1}\ra X^L_{1-23}\cup X^L_{12-3}$}; As $X^R_{1-23}\ra X^L_{1-23}\ra X^R_{23-1}$ we have \sfact{$X^R_{1-23}\ra X^R_{23-1}$}; as $X^R_{13-2}\ra X^L_{23-1}\ra X^R_{1-23}\cup X^R_{23-1}$ we have \sfact{$X^R_{13-2}\ra X^R_{1-23}\cup X^R_{23-1}$}.

Because $X^R_{13-2}\cup X^R_{23-1}\cup X^R_{1-23} \ra X^L_{12-3}$ and $X^R_{123}\ra X^L_{1-23}$, every vertex in $V(R)\setminus Z^R$  has an out-neighbour in $V(L)$. As above, we derive \sfact{$Z^R=\emptyset$}. Similarly, because $X^R_{13-2} \ra X^L_{12-3} \cup X^L_{23-1}$, $X^R_{1-23}\ra X^L_{123}$, and $X^R_{123} \ra X^L_{1-23}$,  every vertex in $V(L)\setminus Z^L$ has in-neighbour in $V(R)$, and so \sfact{$Z^L=\emptyset$}. Next observe that at least one of the sets $X^R_{123},X^L_{123}$ must be empty as $X^R_{123}\ra X^L_{123}\ra X^R_{23-1} \ra X^L_{12-3}\ra X^R_{123}$. If $X^R_{123}=\emptyset$ then $V(R)=X^R_{1-23}\cup X^R_{13-2}\cup X^R_{23-1}$  where each of these sets induces an acyclic subtournament of $R$ and $X^R_{1-23}\rad X^R_{23-1}\rad X^R_{13-2}$ and $X^R_{1-23}\rad X^R_{13-2}$. Thus $R$ is acyclic, contradicting $\inv(R)\geq 2$. So \sfact{$X^R_{123}\neq\emptyset$} and \sfact{$X^L_{123}=\emptyset$}. As above we obtain a contradiction by observing that $L$ is acyclic, contradicting $\inv(L)\geq 2$.
This completes the proof of Subcase 3.1.\\

      \noindent{}{\bf Subcase 3.2} $X^R_{1-23}=\emptyset$.\\
      
      By symmetry, we can assume \fact{$X^R_{2-13}=\emptyset$} and \fact{$X^R_{3-12}\neq\emptyset$}. Hence \fact{$X^L_{1-23}\ra X^L_{3}$} because $X^L_{1-23}\ra X^R_{3-12}\ra X^L_{3}$, and \fact{$X^R_1\ra X^R_{3-12}$} because $X^R_1\ra X^L_{1-23}\ra X^R_{3-12}$. Note that one of $X^L_{13-2},X^R_{13-2}$ is empty since $X^L_{13-2}\ra X^R_{13-2}\ra X^L_{1-23}\ra X^R_{3-12}\ra X^L_{13-2}$. By symmetry we can assume that \fact{$X^L_{13-2}=\emptyset$}. By Claim C, $X^L_2\neq X^L_3$, so \fact{$X^L_{12-3}\neq\emptyset$}. \\

      Suppose first that \asm{$X^R_{123}\neq\emptyset$}. Then \sfact{$X^L_{23-1}=\emptyset$} since $X^L_{23-1}\ra X^R_{123}\ra X^L_{1-23}\ra X^L_{23-1}$. Now, by Claim A, $X^L_3\neq \emptyset$ so \sfact{$X^L_{123}\neq\emptyset$}. 
      Now $X^L_{123} \ra X^R_{13-2}\ra X^L_{1-23}\ra X^L_{123}$, so \sfact{$X^R_{13-2} =\emptyset$}.
      Furthermore, $X^R_{23-1}\ra X^L_{12-3}\ra X^R_{123}\ra X^L_{1-23}\ra X^R_{23-1}$ so \sfact{$X^R_{23-1}=\emptyset$}. 
      Therefore $X^R_1=X^R_2$, a contradiction to Claim C.
      Thus \fact{$X^R_{123}=\emptyset$}.\\

      Next suppose \asm{$X^L_{123}\neq\emptyset$}. Then \sfact{$X^R_{12-3}=\emptyset$} because $X^R_{12-3}\ra X^R_{3-12}\ra X^L_{123}\ra X^R_{12-3}$. By Claim~A, $X^R_1, X^R_2\neq \emptyset$, so  \sfact{$X^R_{13-2}\neq\emptyset$ and \sfact{$X^R_{23-1}\neq\emptyset$}.}    As $X^R_{13-2}\ra X^L_{1-23}\ra X^R_{3-12}\cup X^R_{23-1}$ we have \sfact{$X^R_{13-2}\ra X^R_{3-12}\cup X^R_{23-1}$}. Since $d^+_R(X^R_{13-2})>0$ we have \sfact{$Z^R\neq\emptyset$}. However, there can be no arcs from $Z^R$ to $X^R_3 =V(R)\setminus Z^R$, because $X^R_3\ra X^L_{123}\ra Z^R$.
      This contradicts the fact that $R$ is strong. 
      Thus \fact{$X^L_{123}=\emptyset$}. \\

      By Claim A, $X^L_3\neq \emptyset$, so \fact{$X^L_{23-1}\neq\emptyset$}. Thus   \fact{$X^R_{23-1}=\emptyset$} because $X^R_{23-1}\ra X^L_{12-3}\ra X^R_{3-12}\ra X^L_{23-1}\ra X^R_{23-1}$. 
  By Claim~A,   $X^R_2\neq \emptyset$ so \sfact{$X^R_{12-3}\neq\emptyset$}.  By Claim~C, $X^R_1\neq X^R_2$, so  $X^R_{13-2}\neq\emptyset$. 
As $X^L_{12-3}\ra X^R_{12-3}\ra X^L_{1-23}\cup X^L_{23-1}$, we have \fact{$X^L_{12-3}\ra X^L_{1-23}\cup X^L_{23-1}$}. 
Thus the fact that $d^+_L(X^L_{12-3})>0$ implies that there is an arc $vz$ from $X^L_{12-3}$ to $Z^L$. But then for any $u\in X^R_{13-2}$,  $(u,v,z,u)$ is directed $3$-cycle, a contradiction.

     This completes Subcase 3.2.\\

      \noindent{}{\bf Case 4:} All three  of $X^L_{1-23},X^L_{2-13},X^L_{3-12}$ or all three of   $X^R_{1-23},X^R_{2-13},X^R_{3-12}$ are non-empty.\\

      By symmetry, we can  assume that \fact{$X^L_{1-23},X^L_{2-13},X^L_{3-12}\neq\emptyset$}. There do not exist 
      $i\neq j\in [3]$ such that $X^R_i\setminus X^R_j, X^R_j\setminus X^R_i\neq\emptyset$, for otherwise  $X^L_{i-jk}\ra (X^R_j\setminus X^R_i)\ra X^L_{j-ik}\ra (X^R_i\setminus X^R_j)\ra X^L_{i-jk}$, a contradiction. 
      Hence we may assume by symmetry that \fact{$X^R_2\setminus X^R_1,X^R_3\setminus X^R_1,X^R_2\setminus X^R_3=\emptyset$}.
      This implies that $X^R_2=X^R_{123}$, $X^R_3=X^R_{123}\cup X^R_{13-2}$ and $X^R_1= X^R_{123}\cup X^R_{13-2}\cup X^R_{1-23}$. Moreover,   \fact{$X^R_{1-23},X^R_{123},X^R_{13-2}\neq\emptyset$} by Claim C.
      As $X^R_3\ra X^L_{3-12}\ra X^R_{1-23}$ we have \fact{$X^R_3\ra X^R_{1-23}$}, so since  $d^-_R(X^R_{1-23})>0$ we must have an arc from $Z^R$ to $X^R_{1-23}$ and now $X^R_{1-23}\ra X^L_{1-23}\ra Z^R$ gives a contradiction. This completes Case 4.\\

      \noindent{}{\bf Case 5:}  Exactly two of $X^L_{1-23},X^L_{2-13},X^L_{3-12}$ or two of   $X^R_{1-23},X^R_{2-13},X^R_{3-12}$ are non-empty.\\
      By symmetry we can  assume that \fact{$X^L_{1-23},X^L_{2-13}\neq\emptyset$} and \fact{$X^L_{3-12}=\emptyset$.}\\

      \noindent{}{\bf Subcase 5.1:} $X^R_{1-23},X^R_{2-13},X^R_{3-12}=\emptyset$.\\

    As  $X^L_{1-23}\ra X^R_{23-1} \ra X^L_{2-13}\ra X^R_{13-2} \ra X^L_{1-23}$, one of $X^R_{13-2},X^R_{23-1}$ is empty. 
     By symmetry we may assume that \fact{$X^R_{23-1}=\emptyset$}. 
     By Claim C, $X^R_1 \neq X^R_2$ and $X^R_1 \neq X^R_3$, so \fact{$X^R_{13-2}\neq\emptyset$ and $X^R_{12-3}\neq\emptyset$}.
          Now $V(R)\setminus Z^R=X^R_1\ra X^L_{1-23}\ra Z^R$,  thus there is no arc leaving $Z^R$. As $R$ is strong, we get \fact{$Z^R=\emptyset$}. 
        
  As $X^R_{12-3} \ra X^L_{2-13} \ra X^R_{13-2}$, we have \fact{$X^R_{12-3} \ra X^R_{13-2}$}. Hence as $R$ is strong, necessarily \fact{$X^R_{123}\neq \emptyset$}. If $X^L_{123}\neq\emptyset$, then  $X^R_{123}\ra X^R_{12-3}\cup X^R_{13-2}$ as $X^R_{123}\ra X^L_{123}\ra X^R_{12-3}\cup X^R_{13-2}$. This contradicts the fact that $R$ is strong since $d^+_R(X^R_{123})=0$. Hence \fact{$X^L_{123}=\emptyset$}. By Claim A, $X^L_3\neq \emptyset$, so \fact{$X^L_{13-2}\cup X^L_{23-1}\neq \emptyset$}.

     Since $X^R_{123} \ra X^L_{2-13} \ra X^R_{13-2}$, we have \fact{ $X^R_{123} \ra X^R_{13-2}$}. We also have \fact{$X^R_{12-3}\ra X^R_{123}$} because $X^R_{12-3}\ra X^L_{13-2}\cup X^L_{23-1}\ra X^R_{123}$.
Hence   $V(R)=X^R_{12-3}\cup X^R_{13-2}\cup X^R_{123}$  where each of these sets induces an acyclic subtournament of $R$ and $X^R_{13-2}\rad X^R_{12-3}\rad X^R_{123}\rad X^R_{13-2}$. Thus inverting $X^R_{12-3}\cup X^R_{13-2}$ makes $R$ acyclic, contradicting $\inv(R)\geq 2$.

This completes Subcase 5.1\\

      \noindent{}{\bf Subcase 5.2:} $X^R_{1-23}\neq\emptyset$ and $X^R_{2-13}\cup X^R_{3-12}=\emptyset$.\\

      We first observe that since  $X^L_{2-13}\cup X^L_{23-1}\ra X^R_{1-23}\ra X^L_1$ we can conclude that 
      \fact{$X^L_{2-13}\ra X^L_1$ and $X^L_{23-1}\ra X^L_1$}.
      As $X^R_{23-1}\ra X^L_{2-13}\ra X^L_{1-23}\ra X^R_{23-1}$, we have \fact{$X^R_{23-1}=\emptyset$}. Now $V(R)\setminus Z^R =X^R_1$ and $X^R_1\ra X^L_{1-23} \ra Z^R$. So $V(R)\setminus Z^R \ra Z^R$. Since $R$ is strong, \fact{$Z^R=\emptyset$}.  Now Claims A and C imply that at least two of the sets $X^R_{13-2},X^R_{123},X^R_{12-3}$ are non-empty. This implies that every vertex of $V(L)$ has an in-neighbour in $V(R)$ (as $X^R_{1-23}\ra X^L_1$, $X^R_{13-2}\cup X^R_{12-3}\ra X^L_{23-1}$ and $X^R_2\ra X^L_{2-13}$) so we must have
      \fact{$Z^L=\emptyset$}.

      Suppose first that \asm{$X^R_{12-3}=\emptyset$}. 
      By Claim~A, $X^R_2\neq \emptyset$, so \sfact{$X^R_{123}\neq\emptyset$}. 
      Moreover, by Claim~C, $X^R_2\neq X^R_3$, so \sfact{$X^R_{13-2}\neq\emptyset$}.
       Since $X^L_{12-3}\cup X^L_{13-2}\ra X^R_{123}\ra X^L_{2-13}\ra X^L_{12-3}\cup X^L_{13-2}$ we have \sfact{$X^L_{12-3}\cup X^L_{13-2}=\emptyset$}. If $X^L_{23-1}\neq \emptyset$, then $X^L_{123}=\emptyset$ as $X^L_{23-1}\ra X^L_{123}\ra X^R_{13-2}\ra X^L_{23-1}$ and we have $X^L_{2-13}\ra X^L_{23-1}$ as $X^L_{2-13}\ra X^R_{13-2}\ra X^L_{23-1}$. Now we see that $d^-_L(X^L_{23-1})=0$, a contradiction. Hence \sfact{$X^L_{23-1}=\emptyset$} and \sfact{$X^L_{123}\neq\emptyset$} because $X^L_3\neq \emptyset$ by Claim~A. Moreover \sfact{$X^L_{123}\ra X^L_{1-23}$} because $X^L_{123}\ra X^R_{13-2}\ra X^L_{1-23}$.
     Now $V(L)=X^L_{1-23}\cup X^L_{2-13}\cup X^L_{123}$ where each of these sets induces an acyclic subtournament in $L$ 
     and  $X^L_{1-23}\rad X^L_{123}\rad X^L_{2-13}\rad X^L_{1-23}$. Then inverting the set $X^L_{1-23}\cup X^L_{2-13}$ makes $L$ acyclic, a contradiction to $\inv(L)\geq 2$. Thus \fact{$X^R_{12-3}\neq \emptyset$}.\\

      Note that $X^R_{12-3}\ra X^R_{1-23}\cup X^R_{13-2}$ as $X^R_{12-3}\ra X^L_{2-13}\ra X^R_{1-23}\cup X^R_{13-2}$. Thus \sfact{$X^L_{123}=\emptyset$} because $X^L_{123}\ra X^R_{12-3}\ra X^R_{1-23}\ra X^L_{123}$. Furthermore the fact that $d^+_R(X^R_{12-3})>0$ implies that \fact{$X^R_{123}\neq\emptyset$} and that there is at least one arc from $X^R_{12-3}$ to $X^R_{123}$ in $T$ (and in $R$). We saw before that $X^R_{12-3}\ra X^R_{1-23}$ and by the same reasoning $X^R_{123}\ra X^R_{1-23}$, hence, as $Z^R=\emptyset$ and $d^-_R(X^R_{1-23})>0$, there is at least one arc from $X^R_{1-23}$ to $X^R_{13-2}$. Hence \fact{$X^R_{13-2}\neq\emptyset$} and \sfact{$X^L_{23-1}=\emptyset$} as $X^R_{13-2}\ra X^L_{23-1}\ra X^R_{1-23}$.
        We have \sfact{$X^L_{12-3}=\emptyset$} since $X^L_{12-3}\ra X^R_{123}\ra X^L_{2-13}\ra X^R_{1-23}\ra X^L_{12-3}$. Finally, as $X^L_{2-13}\ra X^R_{1-23}\ra X^L_1$ we have \sfact{$X^L_{2-13}\ra X^L_1$}. But now $d^+_L(X^L_1)=0$ (recall that $Z^L=\emptyset$), a contradiction. This completes Subcase 5.2\\

      \noindent{}{\bf Subcase 5.3:} $X^R_{3-12}\neq\emptyset$ and $X^R_{1-23}\cup X^R_{2-13}=\emptyset$.\\

      As $X^R_{23-1}\ra X^L_{2-13}\ra X^R_{13-2}\ra X^L_{1-23}\ra X^R_{23-1}$ one of the sets $X^R_{13-2},X^R_{23-1}$ must be empty. By symmetry we may assume that \fact{$X^R_{23-1}=\emptyset$}. \\
      
      Suppose first that \asm{$X^R_{12-3}=\emptyset$}. Then, by Claim A, $X^R_2\neq \emptyset$, so \sfact{$X^R_{123}\neq\emptyset$}, and by Claim C, $X^R_1\neq X^R_2$, so \sfact{$X^R_{13-2}\neq\emptyset$}. Now \sfact{$X^L_{123}=\emptyset$} because $X^L_{123}\ra X^R_{13-2}\ra X^L_{1-23}\ra X^R_{3-12}\ra X^L_{123}$. As $X^R_{123}\ra X^L_{2-13}\ra X^R_{13-2}\cup X^R_{3-12}$, we have \sfact{$X^R_{123}\ra X^R_{13-2}\cup X^R_{3-12}$}.
      Next we observe that \sfact{$X^L_{13-2}=\emptyset$} since $X^L_{13-2}\ra X^R_{123}\ra X^R_{3-12}\ra X^L_{13-2}$. Now, as $X^L_3\neq\emptyset$ by Claim C, we have \sfact{$X^L_{23-1}\neq\emptyset$} but that contradicts that  $X^L_{23-1}\ra X^R_{123}\ra X^R_{3-12}\ra X^L_{23-1}$. So we must have \fact{$X^R_{12-3}\neq\emptyset$}.\\

      First observe  that \fact{$X^L_{123}=\emptyset$} as $X^L_{123}\ra X^R_{12-3}\ra X^L_{1-23}\ra X^R_{3-12}\ra X^L_{123}$. As $X^R_1\neq X^R_2$ by Claim~C, we have \fact{$X^R_{13-2}\neq\emptyset$}. Now \fact{$X^L_{13-2}=\emptyset$}
      as $X^L_{13-2}\ra X^R_{13-2}\ra X^L_{1-23}\ra X^R_{3-12}\ra X^L_{13-2}$.
      As $X^L_3\neq\emptyset$ by Claim~A, we have \fact{$X^L_{23-1}\neq\emptyset$.} 
      Since $X^L_{12-3}\ra X^R_{12-3}\ra X^L_{2-13}\ra X^R_{13-2}\ra  X^L_{12-3}$ we have \fact{$X^L_{12-3}=\emptyset$}. 
      As $X^L_{1-23}\ra X^R_{3-12}\ra X^L_{23-1}$, we have \fact{$X^L_{1-23}\ra X^L_{23-1}$}. 
      Moreover $X^L_{2-13}\ra X^R_{13-2}\ra X^L_{23-1}\cup X^L_{1-23}$ implies \fact{$X^L_{2-13}\ra X^L_{23-1}\cup X^L_{1-23}$}. 
      We also have \fact{$Z^L=\emptyset$} since every vertex in $X^L_{1-23}\cup X^L_{23-1}\cup X^L_{2-13}$ has an in-neighbour in $R$, implying that there can be no arc entering $Z^L$.  
      Now $V(L)=X^L_{1-23}\cup X^L_{23-1}\cup X^L_{2-13}$ where each of these sets induces a transitive subtournament in $L$ and  $X^L_{1-23}\rad X^L_{23-1}\rad X^L_{2-13}\rad X^L_{1-23}$. However this implies that
      inverting $X^L_{1-23}\cup X^L_{2-13}$ makes $L$ acyclic, a contradiction to $\inv(L)\geq 2$. This completes the proof of Subcase 5.3.\\

      \noindent{}{\bf Subcase 5.4:} $X^R_{1-23},X^R_{2-13}\neq\emptyset$ and $X^R_{3-12}=\emptyset$.\\

      This case is trivial as $X^L_{1-23}\ra X^R_{2-13}\ra X^L_{2-13}\ra X^R_{1-23}\ra X^L_{1-23}$ contradicts that $T$ is acyclic.\\

      By symmetry the only remaining case to consider is the following.\\
      
      \noindent{}{\bf Subcase 5.5:} $X^R_{1-23},X^R_{3-12}\neq\emptyset$ and $X^R_{2-13}=\emptyset$.\\

      As $X^L_{23-1}\ra X^R_{1-23}\ra X^L_{1-23}\ra X^R_{3-12}\ra X^L_{23-1}$ we have \fact{$X^L_{23-1}=\emptyset$} and as $X^R_{23-1}\ra X^L_{2-13}\ra X^R_{1-23}\ra X^L_{1-23}\ra X^R_{23-1}$ we have \fact{$X^R_{23-1}=\emptyset$}. 
 Note that every vertex in $V(L)$  has an in-neighbour in $V(R)$ (as $X^R_{1-23}\ra X^L_1$ and $X^R_2\ra X^L_{2-13}$) and every vertex in $V(R)$ has an out-neighbour in $V(L)$ (as $X^R_1\ra X^L_{1-23}$ and $X^R_{3-12}\ra X^L_3$). This implies that \fact{$Z^L=\emptyset$} and \fact{$Z^R=\emptyset$}.      
      At least one of $X^L_{13-2}, X^R_{13-2}$ is empty as $X^L_{13-2}\ra X^R_{13-2}\ra X^L_{1-23}\ra X^R_{3-12}\ra X^L_{13-2}$ and at least one of $X^L_{12-3},X^R_{12-3}$ is empty as $X^L_{12-3}\ra X^R_{12-3}\ra X^L_{2-13}\ra X^R_{1-23}\ra X^L_{12-3}$.
     \\

      Suppose first that $X^R_{12-3}=\emptyset=X^R_{13-2}$. Then $X^R_2\neq \emptyset$ by Claim~A, so $X^R_{123}\neq\emptyset$. 
      
    Moreover $X^R_{123}\ra X^R_{1-23}\cup  X^R_{3-12}$ because $X^R_{123}\ra X^L_{2-13}\ra X^R_{1-23}\cup  X^R_{3-12}$. This implies that $d^+_R(X^R_{123})=0$, a contradiction.\\

      Suppose next that $X^L_{12-3}=\emptyset= X^L_{13-2}$. Then $X^L_3\neq \emptyset$ by Claim~A, so $X^L_{123}\neq\emptyset$.
Moreover $X^L_{1-23}\cup X^L_{2-13}\ra X^L_{123}$ as $X^L_{1-23}\cup X^L_{2-13}\ra X^R_{3-12}\ra X^L_{123}$.
 This implies that $d^-_L(X^L_{123})=0$, a contradiction.\\

      Now assume that $X^R_{12-3}=\emptyset=X^L_{13-2}$ and $X^R_{13-2}\neq\emptyset\neq X^L_{12-3}$. Then $X^L_{123}\neq\emptyset$ as $X^L_3\neq\emptyset$ by Claim A and now we get the contradiction
      $X^L_{123}\ra X^R_{13-2}\ra X^L_{1-23}\ra X^R_{3-12}\ra X^L_{123}$.\\

      The final case is $X^R_{12-3}\neq\emptyset\neq X^L_{13-2}$ and $X^R_{13-2}=\emptyset= X^L_{12-3}$.
      We first observe that $X^R_{123}=\emptyset$ as $X^R_{123}\ra X^L_{1-23}\ra X^R_{3-12}\ra X^L_{13-2}\ra X^R_{123}$. As $X^R_{12-3}\ra X^L_{2-13}\ra X^R_{1-23}$ we have $X^R_{12-3}\ra X^R_{1-23}$ and as $X^R_{1-23}\ra X^L_{1-23}\ra X^R_{3-12}$ we have $X^R_{1-23}\ra X^R_{3-12}$. This implies that $d^-_R(X^R_{1-23})=0$, a contradiction.  This completes the proof of Subcase 5.5 and the proof of the theorem.

  \end{proof}

  \begin{corollary}\label{cor:2+2=4strong}
  Let $L$ and $R$ be strong oriented graphs such that $\inv{}(L),\inv{}(R) = 2$. Then $\inv{}(L\rightarrow{}R) = 4$.
  \end{corollary}

 \section{Inversion number of augmentations of oriented graphs}\label{sec:augment}

\begin{lemma}\label{lem:auginv=1}
 Let $D$ be an oriented graph  with $\inv(D)=1$.
 Then $\inv(\sigma(z,D)) = 2$ for every $z\in V(D)$.
\end{lemma}
\begin{proof}
Recall that $\inv(\sigma(z,D)) \leq \inv(D)+1 = 2$ for every vertex $z\in V(D)$.

Suppose for a contradiction that there is a vertex $z$ of $D$ such that $\inv(\sigma(z,D)) = 1$. Let $X$ be a set whose inversion in $\sigma(z,D)$ results in an acyclic digraph $D^*$.  

As $D$ has inversion number $1$ it has a directed cycle $C$. The set $X$ contains an arc $uu^+$ of $C$, for otherwise $C$ would be a directed cycle in $D^*$. Moreover, $X$ does not contain all vertices of $C$, for otherwise the inversion of $X$ transforms $C$ in the directed cycle in the opposite direction. Hence, without loss of generality, we may assume that $u^-$, the in-neighbour of $u$ in $C$ is not in $X$.

Note also that $C'=(z,y,x,z)$ is a directed cycle in $\sigma(z,D)$ so $X$ must contain exactly two vertices of $C'$.
 In particular, there is a vertex, say $w$, in $\{x,y\}\cap X$.
 \begin{itemize}
\item  If $z\notin \{u^-, u\}$, then $(w, u^-, u, w)$ is a directed $3$-cycle, a contradiction.
 
\item If $z=u$, then either $X\cap V(C')=\{x,z\}$ and $(z,x,u^-,z)$ is a directed $3$-cycle in $D^*$, or $X\cap V(C')=\{y,z\}$ and $(x,u^+,y,x)$ is a directed $3$-cycle in $D^*$, a contradiction.

\item If $z=u^-$, then $X\cap V(C')=\{x,y\}$ and $(z,u,x,z)$  is a directed $3$-cycle in $D^*$, a contradiction.
\end{itemize}
  \end{proof}

 Recall that $\sigma_i(z, D)$ denotes the $z$-augmentation of $D$ on which the vertices added are denoted by $x_i$ and $y_i$.
 \begin{theorem}\label{thm:2augment}
	Let $D$ be an oriented graph with $\inv(D)=1$ and let $H=\sigma_1(x_2,\sigma_2(z, D))$. Then, $\inv(H)=3$.
 \end{theorem}
 \begin{proof}
 By Lemma~\ref{lem:auginv=1}, $\inv(\sigma_2(z, D)) =2$. In addition, $\sigma_2(z, D)$ is a subdigraph of $H$, 
 so by Proposition~\ref{prop:monotone}, $\inv(H)\geq 2$. Moreover, $\inv(H) \leq \inv(\sigma_2(z, D)) +1 =3$.
 
 \medskip
 
 Assume for a contradiction that $\inv(H) = 2$.
 Let $(X_1, X_2)$ be a decycling family of $H$. For $i\in [2]$, let $H_i= \Inv(H; X_i)$. Note that $\inv(H_i)\leq 1$ for  $i\in [2]$, because $(X_1,X_2)$ is a decycling family.
 
 Then $(X_1\setminus \{y_2\}, X_2\setminus \{y_2\})$ is a decycling family of $H-y_2$. 
 But $H-y_2$ is isomorphic to $ \vec{C_3}\rightarrow D$ with $(y_1, x_1, x_2, y_1)$ dominating $D$.
 Thus, by Proposition \ref{prop:1>1}, $\inv(H-y_2)\geq 2$ and furthermore, by Theorem~\ref{thm:1>1}, we may assume that $X_1 \subseteq\{x_1, y_1, x_2, y_2\}$. Observe that $X_1\cap |\{x_1, y_1, x_2\}|=2$, for otherwise $H_1-y_2 = H-y_2$ and $\inv(H_1)\geq \inv(H_1-y_2)\geq 2$ by Proposition~\ref{prop:monotone}.
Hence, there is a vertex $v\in \{y_1, x_1, x_2\}$ such that $y_2v\in A(H_1)$.
This implies that  $H_1\langle \{v, y_2\}\cup V(D)\rangle$ is a $z$-augmentation of $D$. Therefore, by Lemma~\ref{lem:auginv=1}, we must have $\inv(H_1\langle \{v, y_2\}\cup V(D)\rangle)=2$, and so by Proposition~\ref{prop:monotone},  $\inv(H_1)\geq 2$, a contradiction.

Thus, we have shown that $\inv(H)=3$. 
\end{proof}

Recall that $Q_n$ is the tournament we obtain from the transitive tournament on $n$ vertices by reversing the arcs of the unique hamiltonian path $(v_1, \ldots{},v_n)$. Hence $Q_7$ is the oriented graph $\sigma{}_1(v_3,\sigma{}_2(v_5,\vec{C_3})$, where $\vec{C_3}$ is the  directed $3$-cycle $(v_5,v_7,v_6,v_5)$, $x_2=v_3$, $y_2=v_4$, $x_1=v_1$ and $y_1=v_2$. Thus Theorem~\ref{thm:2augment} yields the following.
 
 \begin{corollary}
   \label{cor:invQ7}
 $\inv(Q_7)=3$.
\end{corollary}

\section{Inversion number of intercyclic oriented graphs}\label{sec:intercyclic}

A digraph $D$ is {\bf intercyclic} if $\nu(D) = 1$. 
The aim of this subsection is to prove the following theorem.

\begin{theorem}\label{thm:intercyclic}
If $D$ is an intercyclic oriented graph, then $\inv(D)\leq 4$.
\end{theorem}

In order to prove this theorem, we need some preliminaries.

Let $D$ be an oriented graph.
An arc $uv$ is {\bf weak} in $D$ if $\min \{d^+(u), d^-(v)\}=1$.
An arc is {\bf  contractable} in $D$ if it is weak and in no directed $3$-cycle.
If $a$ is a contractable arc, then let $D/a$ is the digraph obtained
by contracting the arc $a$ and $\tilde{D}/a$ be the oriented graph obtained from $D$ by removing one arc from every  pair of parallel arcs created in $D/a$.

\begin{lemma}\label{lem:contraction}
Let $D$ be a strong oriented graph and let $a$ be a contractable arc in $D$.
Then $D/a$ is a strong intercyclic oriented graph and $\inv(\tilde{D}/a) \geq \inv(D)$.
\end{lemma} 
\begin{proof}
McCuaig proved that $D/a$ is strong and intercyclic.
Let us prove that  $\inv(D) \leq \inv(\tilde{D}/a)$.
Observe that $\inv(\tilde{D}/a) = \inv(D/a)$.

Set $a=uv$, and let $w$ be the vertex corresponding to both $u$ and $v$ in $D/a$.
Let $(X'_1, \dots , X'_p)$ be a decycling family of $D'=\tilde{D}/a$ that result in an acyclic oriented graph $R'$.
For $i\in [p]$, let $X_i = X'_i$ if $w\notin X'_i$ and $X_i = (X'_i\setminus \{w\}) \cup \{u,v\}$ if $w\in X'_i$.
Let $a^*=uv $ if $w$ is in an even number of $X'_i$ and $a^*=vu$ otherwise, and let $R = \Inv(D; (X_1, \dots , X_p))$.
One easily shows that $R = R'/a^*$. Therefore $R$ is acyclic since the contraction of an arc transforms a directed cycle into a directed cycle. \end{proof}

\begin{lemma}\label{lem:non-contract}
Let $D$ be an intercyclic oriented graph.
If there is a non-contractable weak arc, then $\inv(D) \leq 4$.
\end{lemma}

\begin{proof}
Let $uv$ be a non-contractable weak arc. By directional duality, we may assume that $d^-(v)=1$.
Since $uv$ is non-contractable, $uv$ is in a directed $3$-cycle $(u,v,w,u)$.
Since $D$ is intercyclic, we have $D\setminus \{u,v,w\}$ is acyclic.
Consequently, $\{w,u\}$ is a cycle transversal of $D$, because every directed cycle containing $v$ also contains $u$.
Hence, by Theorem~\ref{thm:bound-fvs}, $\inv(D) \leq 2\tau(D) \leq 4$.
\end{proof}

The description below follows \cite{bangC31}.
A digraph $D$ is {\bf in reduced form} if it is strong, and it has no weak arc, that is $\min \{\delta^-(D), \delta^+(D)\} \geq 2$.

\begin{figure}[hbtp]
\begin{center}
\tikzstyle{vertex}=[circle,draw, minimum size=14pt, scale=0.6, inner sep=0.5pt]
\begin{tikzpicture}[scale=0.7]
  \node(y1) at (0,2) [vertex] {$y_1$};
\node(y2) at (1,4) [vertex] {$y_2$};
\node(y3) at (3,4) [vertex] {$y_3$};
\node(y4) at (4,2) [vertex] {$y_4$};
\node(y5) at (3,0) [vertex] {$y_5$};
\node(y6) at (1,0) [vertex] {$y_6$};
\node(y) at (2,2) [vertex] {$y$};

\node () at (2,-1) {(a)};

\draw[->, line width=0.03cm] (y1) to (y2);
\draw[->, line width=0.03cm] (y2) to (y3);
\draw[->, line width=0.03cm] (y3) to (y4);
\draw[->, line width=0.03cm] (y4) to (y5);
\draw[->, line width=0.03cm] (y5) to (y6);
\draw[->, line width=0.03cm] (y6) to (y1);
\draw[->, line width=0.03cm] (y1) to (y);
\draw[->, line width=0.03cm] (y) to (y2);
\draw[->, line width=0.03cm] (y3) to (y);
\draw[->, line width=0.03cm] (y) to (y4);
\draw[->, line width=0.03cm] (y5) to (y);
\draw[->, line width=0.03cm] (y) to (y6);

\draw[->, line width=0.03cm] (y4) to [out=210, in=-30] (y1);
\draw[->, line width=0.03cm] (y2) to [out=-90,in=150] (y5);
\draw[->, line width=0.03cm] (y6) to [out=30,in=-90] (y3);

\node(y1a) at (6,2) [vertex] {$y_1$};
\node(y2a) at (7,4) [vertex] {$y_2$};
\node(y3a) at (9,4) [vertex] {$y_3$};
\node(y4a) at (10,2) [vertex] {$y_4$};
\node(y5a) at (9,0) [vertex] {$y_5$};
\node(y6a) at (7,0) [vertex] {$y_6$};
\node(ya) at (8,2) [vertex] {$y$};

\node () at (8,-1) {(b)};

\draw[->, line width=0.03cm] (y1a) to (y2a);
\draw[->, line width=0.03cm] (y2a) to (y3a);
\draw[->, line width=0.03cm] (y3a) to (y4a);
\draw[->, line width=0.03cm] (y4a) to (y5a);
\draw[->, line width=0.03cm] (y5a) to (y6a);
\draw[->, line width=0.03cm] (y6a) to (y1a);
\draw[->, line width=0.03cm] (y1a) to (ya);
\draw[->, line width=0.03cm] (y2a) to (ya);
\draw[->, line width=0.03cm] (y3a) to (ya);
\draw[->, line width=0.03cm] (y4a) to (ya);
\draw[->, line width=0.03cm] (y5a) to (ya);
\draw[->, line width=0.03cm] (y6a) to (ya);

\draw[->, line width=0.03cm] (y4a) to [out=210, in=-30] (y1a);
\draw[->, line width=0.03cm] (y2a) to [out=-90,in=150] (y5a);
\draw[->, line width=0.03cm] (y6a) to [out=30,in=-90] (y3a);


\node(y1b) at (12,2) [vertex] {$y_1$};
\node(y2b) at (13,4) [vertex] {$y_2$};
\node(y3b) at (15,4) [vertex] {$y_3$};
\node(y4b) at (16,2) [vertex] {$y_4$};
\node(y5b) at (15,0) [vertex] {$y_5$};
\node(y6b) at (13,0) [vertex] {$y_6$};
\node(yb) at (14,2) [vertex] {$y$};

\node () at (14,-1) {(c)};

\draw[->, line width=0.03cm] (y1b) to (y2b);
\draw[->, line width=0.03cm] (y3b) to (y2b);
\draw[->, line width=0.03cm] (y3b) to (y4b);
\draw[->, line width=0.03cm] (y4b) to (y5b);
\draw[->, line width=0.03cm] (y6b) to (y5b);
\draw[->, line width=0.03cm] (y6b) to (y1b);
\draw[->, line width=0.03cm] (y1b) to (yb);
\draw[->, line width=0.03cm] (y2b) to (yb);
\draw[->, line width=0.03cm] (y3b) to (yb);
\draw[->, line width=0.03cm] (y4b) to (yb);
\draw[->, line width=0.03cm] (y5b) to (yb);
\draw[->, line width=0.03cm] (y6b) to (yb);

\draw[->, line width=0.03cm] (y4b) to [out=210, in=-30] (y1b);
\draw[->, line width=0.03cm] (y5b) to [out=150,in=-90] (y2b);
\draw[->, line width=0.03cm] (y3b) to [out=-90,in=30] (y6b);

\end{tikzpicture}
\end{center}
\caption{(a): the digraph $D_7$; (b): the digraph $D'_7$ obtained from $D_7$ by  inverting the set $\{y,y_2,y_4,y_6\}$; (c): the acyclic digraph  $D''_7$ obtained from $D'_7$ by inverting the set $\{y_2,y_3,y_5,y_6\}$.}\label{fig:D7}
\end{figure}
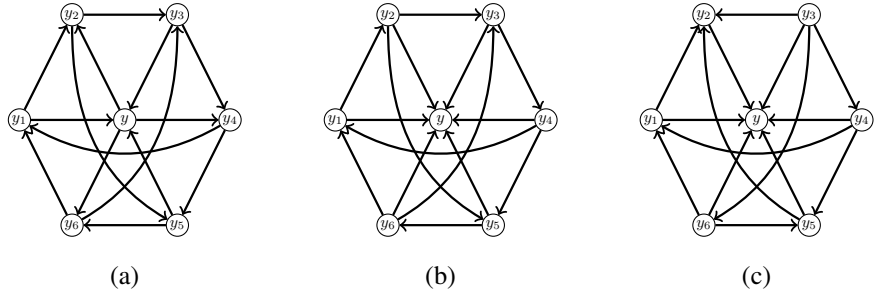

Intercyclic digraphs in reduced form were characterized by Mc Cuaig~\cite{McCuaig91}.
In order to restate his result, we need some definitions.
Let ${\cal P}(x_1,\dots,x_s; y_1,\dots,y_t)$ be the class of acyclic digraphs $D$
such that $x_1,\dots,x_s$, $s \geq 2$, are the sources of $D$, $y_1,\dots,y_t$, $t \geq 2$,
are the sinks of $D$, every vertex which is neither a source nor a sink has in- and out-degree
at least $2$, and, for $1 \leq i<j \leq s$ and $1 \leq k<\ell \leq t$, every $(x_i,y_\ell)$-path intersects
every $(x_j,y_k)$-path. 
By a theorem of Metzlar \cite{Metzlar1989}, such a digraph can be embedded in
a disk such that $x_1,x_2,\dots,x_s,y_t,y_{t-1},\dots,y_1$ occur, in this cyclic order, on its
boundary.
Let ${\cal T}$ be the class of digraphs with minimum in- and out-degree at least $2$
which can be obtained from a digraph in ${\cal P}(x^+,y^+;x^-,y^-)$ by identifying $x^+=x^-$
and $y^+=y^-$.
Let $D_7$ be the digraph from Figure \ref{fig:D7}(a).

\begin{figure}[hbtp]
  \begin{center} \includegraphics[scale=0.6]{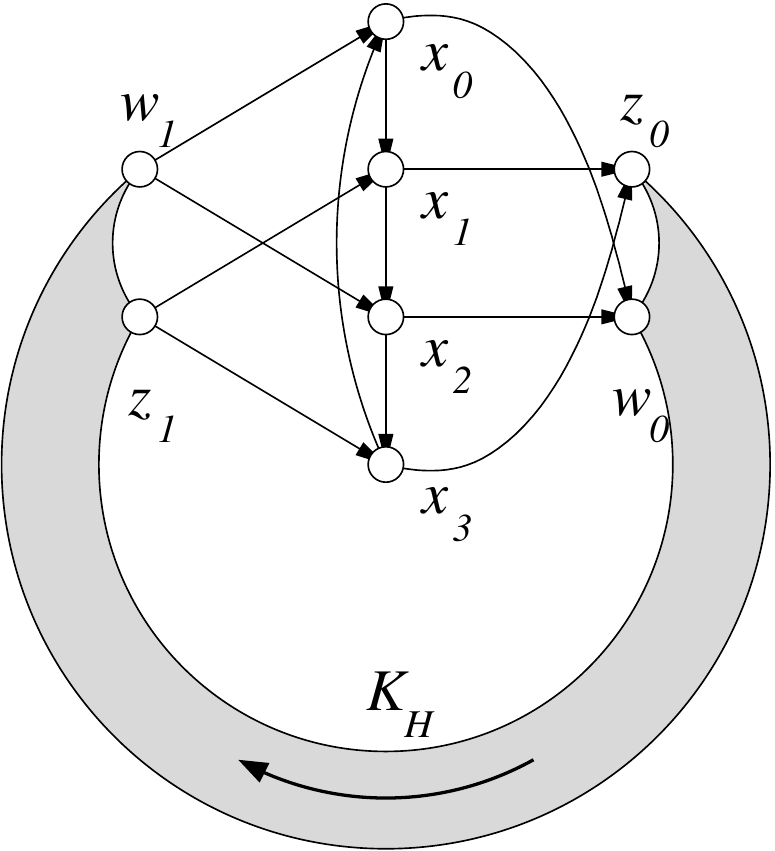} \end{center}
  \caption{\label{F2} The digraphs from ${\cal K}$.
  The arrow in the grey area symbolizing the acyclic (plane) digraph $K_H$ indicates that 
  $z_0,w_0$ are its sources and $z_1,w_1$ are its sinks. (This figure is a courtesy of \cite{bangC31}).}
\end{figure}

\begin{figure}[hbtp]
  \begin{center}
    \includegraphics[scale=0.6]{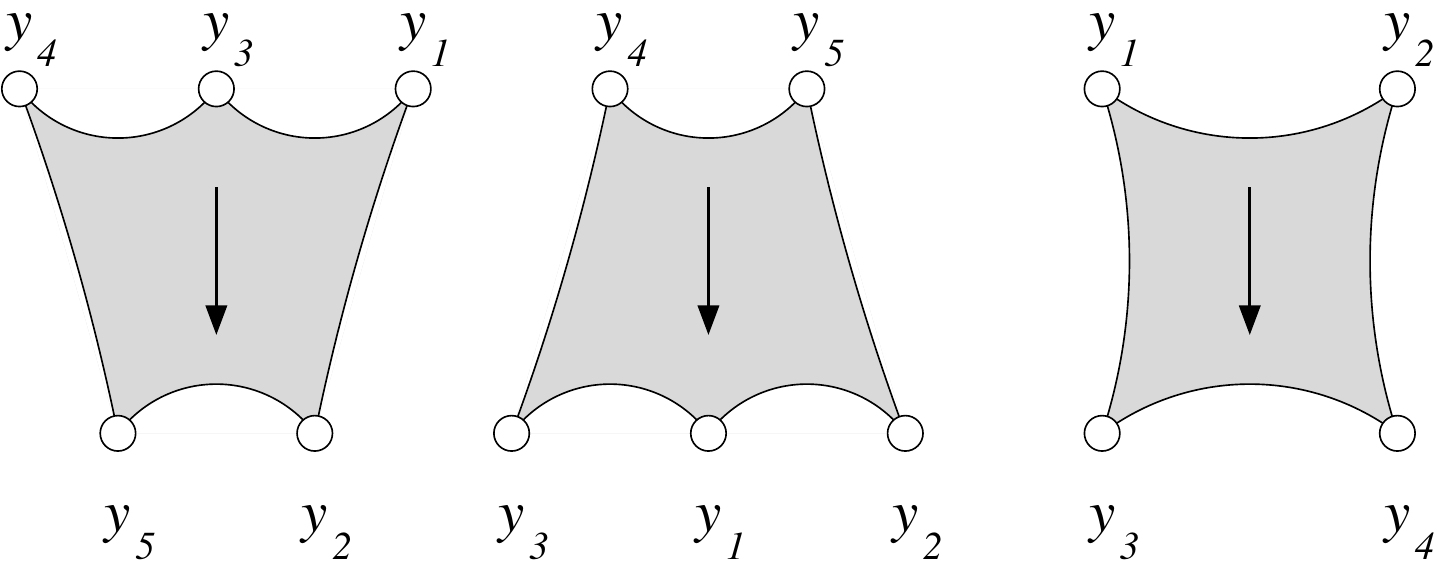} \\[2mm]
    \hspace*{2mm} $H_\alpha$ \hspace{23mm} $H_\beta$ \hspace{25mm} $H_\gamma$
  \end{center}
  \caption{\label{F3} The digraphs from ${\cal H}$. (This figure is a courtesy of \cite{bangC31}).}
\end{figure}

Let ${\cal K}$ be the class of digraphs $D$ with $\tau{}(D) \geq 3$ and
$\delta^0(D) \geq 2$ (Recall that $\delta^0(D) = \min \{\delta^+(D), \delta^-(D)\}$.) which can be obtained from a digraph
$K_H$ from ${\cal P}(w_0,z_0;z_1,w_1)$ by adding at most one arc
connecting $w_0,z_0$, adding at most one arc connecting $w_1,z_1$,
adding a directed $4$-cycle $(x_0, x_1, x_2, x_3, x_0)$ disjoint from $K_H$ and adding eight
single arcs $w_1x_0$, $w_1x_2$, $z_1x_1$, $z_1x_3$, $x_0w_0$, $x_2w_0$, $x_1z_0$ ,$x_3z_0$
(see Figure \ref{F2}).
Let ${\cal H}$ be the class of digraphs $D$ with $\tau{}(D) \geq 3$ and
$\delta^0(D) \geq 2$ such that $D$ is the union of three
arc-disjoint digraphs $H_\alpha \in {\cal P}(y_4,y_3,y_1;y_5,y_2)$,
$H_\beta \in {\cal P}(y_4,y_5;y_3,y_1,y_2)$, and $H_\gamma \in {\cal P}(y_1,y_2;y_3,y_4)$,
where $y_1,y_2,y_3,y_4,y_5$ are the only vertices in $D$ occurring in more than
one of $H_\alpha,H_\beta,H_\gamma$ (see Figure \ref{F3}).

\begin{theorem}[McCuaig~\cite{McCuaig91}]
  \label{T1}
 
   The class of intercyclic digraphs  in reduced form is
   ${\cal T} \cup \{D_7\} \cup {\cal K} \cup {\cal H}$.
\end{theorem}

Using this characterization we can now prove the following. 
\begin{corollary}\label{cor:reduced}
If $D$ is an intercyclic oriented graph in reduced form, then $\inv(D)\leq 4$.
\end{corollary}
\begin{proof}
Let $D$ be an intercyclic oriented graph in reduced form.
By Theorem~\ref{T1}, it is in ${\cal T} \cup \{D_7\} \cup {\cal K} \cup {\cal H}$.

\medskip

If $D\in {\cal T}$, then it is obtained from a digraph $D'$ in ${\cal P}(x^+,y^+;x^-,y^-)$ by identifying $x^+=x^-$
and $y^+=y^-$. Thus $D-\{x^+, y^+\}= D'-\{x^+,y^+, x^-, y^-\}$ is acyclic. Hence $\tau(D)\leq 2$, and  so by Theorem~\ref{thm:bound-fvs}, $\inv(D)\leq 4$.
\medskip

If $D=D_7$, then  inverting $X_1= \{y,y_2,y_4,y_6\}$ so that $y$ becomes a sink and then inverting $\{y_2,y_3,y_5,y_6\}$, we obtain an acyclic digraph with acyclic ordering $(y_3,y_6,y_4,y_5,y_1,y_2,y)$. (See Figure~\ref{fig:D7}.) Hence $\inv(D_7)\leq 2$.

If $D \in {\cal K}$, then inverting $\{x_0,x_3\}$ and $\{x_0,x_1,x_2,x_3,w_1,z_1\}$, we convert $D$ to an acyclic digraph with acyclic ordering $(x_3,x_2,x_1,x_0, v_1, \dots , v_p)$ where  $(v_1, \dots , v_p)$ is an acyclic ordering of $K_H$.

\medskip

If  $D\in {\cal H}$, then consider $D' = \Inv(D, V(H_{\gamma}))$.
The oriented graph $D'$ is the union of   $H_\alpha \in {\cal P}(y_4,y_3,y_1;y_5,y_2)$,
$H_\beta \in {\cal P}(y_4,y_5;y_3,y_1,y_2)$, and $\stackrel{\leftarrow}{H}_\gamma$, the converse of $H_\gamma$.
As  $H_\gamma\in {\cal P}(y_1,y_2;y_3,y_4)$, we have $\stackrel{\leftarrow}{H}_\gamma \in  {\cal P}(y_4,y_3;y_2,y_1)$.
Set $Y=\{y_1, y_2, y_3, y_4, y_5\}$.

We claim that every directed cycle $C'$ of $D'$ contains $y_5$.
Since $D'-Y$ is acyclic, $C'$ is the concatenation of directed paths $P_1, P_2, \dots , P_q$ with both end-vertices in $Y$ and no internal vertex in $Y$. 
Now let $C$ be the directed cycle obtained from $C'$ by replacing each $P_i$ by an arc from its initial vertex to its terminal vertex.
Clearly, $C$ contains $y_5$ if and only if $C'$ does.
But $C$ is a directed cycle in $J$ the digraph with vertex set $Y$ in which $\{y_4, y_3, y_1\}\ra \{y_5,y_2\}$, $\{y_4, y_5\} \ra\{y_3,y_1, y_2\}$, and $\{y_4, y_3\} \ra \{y_1,y_2\}$.
One easily checks that $J-v_5$ is acyclic with acyclic ordering $(y_4,y_3, y_1, y_2)$, so $C$ contains $y_5$ and so does $C'$.

Consequently, $\{y_5\}$ is a cycle transversal of $D'$. Hence, by Theorem~\ref{thm:bound-fvs}, we have $\inv(D')\leq 2\tau(D') \leq 2$.
As $D'$ is obtained from $D$ by inverting one set, we get $\inv(D)\leq 3$.
\end{proof}

 We can now prove Theorem~\ref{thm:intercyclic}.
 
 \begin{proof}
 By induction on the number of vertices of $D$, the result holding trivially if $|V(D)|=3$, that is $D=\vec{C_3}$.

Assume now that $|V(D)|>3$.

 If $D$ is not strong, then it has a unique non-trivial strong component $C$ and any decycling family of $C$ is a decycling family of $D$, so $\inv(C)=\inv(D)$.
 By the induction hypothesis, $\inv(C) \leq 4$, so $\inv(D)\leq 4$.
Henceforth, we may assume that $D$ is strong. 

Assume now that $D$ has a weak arc $a$.
If $a$ is non-contractable, then $\inv(D)\leq 4$ by Lemma~\ref{lem:non-contract}.
If $a$ is contractable, then consider $\tilde{D}/a$. As observed by McCuaig~\cite{McCuaig91}, $D/a$ is also intercyclic.
So by Lemma~\ref{lem:contraction} and the induction hypothesis, $\inv(D)\leq \inv(D/a) \leq 4$.
Henceforth, we may assume that $D$ has no weak arc.

Thus $D$ is in a reduced form and by Corollary~\ref{cor:reduced}, $\inv(D)\leq 4$.
 \end{proof}

\section{Complexity results}\label{sec:complexity}

 \subsection{NP-hardness of {\sc $1$-Inversion} and {\sc $2$-Inversion}}\label{subsec:NP}

 \begin{theorem}\label{thm:NP1}
 {\sc $1$-Inversion} is NP-complete even when restricted to strong oriented graphs.
 \end{theorem}

 In order to prove this theorem, we need some preliminaries.
 
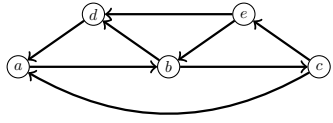
\begin{figure}[hbtp]
\begin{center}
\tikzstyle{vertex}=[circle,draw, minimum size=14pt, scale=0.6, inner sep=0.5pt]

\begin{tikzpicture}[scale=1]
 \node (a) at (0,0) [vertex] {$a$};
  \node (b) at (2,0) [vertex] {$b$};
  \node (c) at (4,0) [vertex] {$c$};
 \node (d) at (1,0.7) [vertex] {$d$};
  \node (e) at (3,0.7) [vertex] {$e$};

  \draw [->, line width=0.03cm] (a) to (b);
\draw [->, line width=0.03cm] (b) to (c);
 \draw [->, line width=0.03cm] (c) to (e);
\draw [->, line width=0.03cm] (e) to (b);
 \draw [->, line width=0.03cm] (d) to (a);
\draw [->, line width=0.03cm] (b) to (d);
\draw [->, line width=0.03cm] (e) to (d);

\draw [->, line width=0.03cm] (c) to [in=-30, out=-150] (a);
\end{tikzpicture}
\caption{The oriented graph $J$}\label{fig:J}
\end{center}
\end{figure}

 Let $J$ be the oriented graph depicted in Figure~\ref{fig:J}.
 
\begin{lemma}\label{lem:J}
The only sets whose inversion can make $J$ acyclic are $\{a,b,e\}$ and $\{b,c,d\}$. 
\end{lemma}
 \begin{proof}
 Assume that an inversion on $X$ makes $J$ acyclic.
 Then $X$ must contain exactly two vertices of each of the directed $3$-cycles $(a,b,c,a)$, $(a,b,d,a)$, and $(e,b,c,e)$, and cannot be $\{a,c,d,e\}$ for otherwise $(e,b,d,e)$ is a directed cycle in the resulting oriented graph.
 Hence $X$ must be either $\{a,b,e\}$ or $\{b,c,d\}$. 
 One can easily check that an inversion on any of these two sets makes $J$ acyclic.
 \end{proof}
 
 \begin{proof}[Proof of Theorem~\ref{thm:NP1}]
 Reduction from {\sc Monotone 1-in-3 SAT} which is well-known to be NP-complete.
 
 Let $\Phi$ be a monotone 3-SAT formula with variables $x_1,\ldots,x_n$ and clauses $C_1,\ldots,C_m$.
 Let $D$ be the oriented graph constructed as follows.
 For every $i\in [n]$, let us construct a variable digraph $K_i$ as follows:
 for every $j\in [m]$, create a copy $J_i^j$ of $J$, and then identify all the vertices $c_i^j$ into one vertex $c_i$ as depicted in Figure \ref{fig:Ki}.
 Then, for every clause $C_j=x_{i_1}\vee x_{i_2} \vee x_{i_3}$, we add the arcs of the directed $3$-cycle $D_j=(a_{i_1}^j, a_{i_2}^j, a_{i_3}^j)$.
 
 \begin{figure}[hbtp]
  \begin{center}
    \tikzstyle{vertex}=[circle,draw, minimum size=14pt, scale=0.6, inner sep=0.5pt]
    
    \begin{tikzpicture}[scale=1]
      \node (a) at (0,0)  {};
      \node (b) at (2,0)  {};
      \node (d) at (1,0.7) {};
      \node (e) at (3,0.7) {};
      \def\xa{0}; \def\ya{0};
      \def\xc{3.5}; \def\yc{0};
      \def\xd{1}; \def\yd{0.7};
      \def\xe{2.5}; \def\ye{0.7};
      
      \draw (\xa,\ya) -- (\xd,\yd) -- (\xe,\ye) -- (\xc,\yc) -- (\xa,\ya);
      \node (l) at (1.75, 0.35) {$J_i^2$};
      \begin{scope}[rotate around ={-45:(\xc, \yc)}]
        \draw (\xa,\ya) -- (\xd,\yd) -- (\xe,\ye) -- (\xc,\yc) -- (\xa,\ya);
        \node (l) at (1.75, 0.35) {$J_i^1$};
      \end{scope}
      
      \begin{scope}[rotate around ={45:(\xc, \yc)}]
        \draw (\xa,\ya) -- (\xd,\yd) -- (\xe,\ye) -- (\xc,\yc) -- (\xa,\ya);
        \node (l) at (1.75, 0.35) {$J_i^3$};
      \end{scope}

      \begin{scope}[rotate around ={75:(\xc, \yc)}]
        \node (l) at (1.75, 0.35) {$\ddots$};
      \end{scope}
      
      \begin{scope}[rotate around ={100:(\xc, \yc)}]
        \draw (\xa,\ya) -- (\xd,\yd) -- (\xe,\ye) -- (\xc,\yc) -- (\xa,\ya);
        \node (l) at (1.75, 0.35) {$J_i^m$};
      \end{scope}
      
      \node[vertex, fill=white] (c) at (\xc,\yc)  {\Large $c_i$};
      
    \end{tikzpicture}
    \caption{The variable gadget $K_i$}\label{fig:Ki}
  \end{center}
  \end{figure}
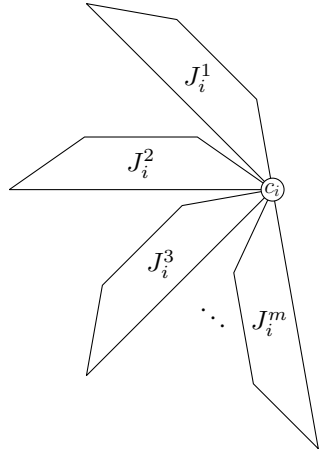

 Observe that $D$ is strong.
 We shall prove that $\inv(D)=1$ if and only if $\Phi$ admits a 1-in-3-SAT assignment.
 
 \medskip
 
 Assume first that $\inv(D)=1$.
 Let $X$ be a set whose inversion makes $D$ acyclic.
 By Lemma~\ref{lem:J}, and the vertices $c_i^j$ are identified in $c_i$, for every $i\in [n]$, either $X\cap V(K_i) = \bigcup_{j=1}^m \{a_i^j, b_i^j, e_i^j\}$ or 
 $X\cap V(K_i) = \bigcup_{j=1}^m \{b_i^j, c_i, d_i^j\}$.
 Let $\varphi$ be the truth assignment defined by $\varphi(x_i) = true$ if $X\cap V(K_i) = \bigcup_{j=1}^m \{b_i^j, c_i, d_i^j\}$, and $\varphi(x_i) = false$ if $X\cap V(K_i) = \bigcup_{j=1}^m \{a_i^j, b_i^j, e_i^j\}$.
 
 Consider a clause $C_j= x_{i_1}\vee x_{i_2} \vee x_{i_3}$. 
 Because $D_j$ is a directed $3$-cycle, $X$ contains exactly two vertices in $V(D_j)$. Let $\ell_1$ and $\ell_2$ be the two indices of $\{i_1, i_2, i_3\}$ such that $a_{\ell_1}^j$ and $a_{\ell_2}^j$ are in $X$ and $\ell_3$ be the third one.
By our definition of $\varphi$, we have $\varphi(x_{\ell_1}) = \varphi(x_{\ell_2}) = false$ and $\varphi(x_{\ell_3}) = true$.
Therefore, $\varphi$ is a 1-in-3 SAT assignment.

 \medskip
 
 Assume now that $\Phi$ admits a 1-in-3 SAT assignment $\varphi$.
For every $i\in [n]$, let $X_i=  \bigcup_{j=1}^m \{b_i^j, c_i, d_i^j\}$ if  $\varphi(x_i) = true$ and $X_i= \bigcup_{j=1}^m \{a_i^j, b_i^j, e_i^j\}$ if $\varphi(x_i) = false$, and set $X=\bigcup_{i=1}^n X_i$.

Let $D'$ be the graph obtained from $D$ by the inversion on $X$.
We shall prove that $D$ is acyclic, which implies $\inv(D)=1$.

Assume for a contradiction that $D'$ contains a directed cycle $C$.
By Lemma~\ref{lem:J}, there is no directed cycle in any variable gadget $K_i$, so $C$ must contain an arc with both ends in $V(D_j)$ for some $j$.
Let $C_j = x_{i_1}\vee x_{i_2} \vee x_{i_3}$.
Now since $\varphi$ is a 1-in-3-SAT assignment, w.l.o.g., we may assume that $\varphi(x_{i_1}) = \varphi(x_{i_2}) = false$ and $\varphi(x_{i_3}) = true$.
Hence in $D'$, $a^j_{i_2}\ra a^j_{i_1}$, $a^j_{i_2}\ra a^j_{i_3}$ and $a^j_{i_3}\ra a^j_{i_1}$.
Moreover, in $D'\langle V(J^j_{i_1})\rangle$, $a^j_{i_1}$ is a sink, so $a^j_{i_1}$ is a sink in $D'$.
Therefore $C$ does not goes through $a^j_{i_1}$, and thus $C$ contains the arc $a^j_{i_2}a^j_{i_3}$, and then enters
$J_{i_3}^j$. But in $D'\langle V(J_{i_3}^j)\rangle$, $a^j_{i_3}$ has a unique out-neighbour, namely $b^j_{i_3}$, which is a sink. This is a contradiction.
\end{proof}

  \begin{corollary}
    \label{cor:NP2}
{\sc $2$-Inversion} is NP-complete.
  \end{corollary}
  \begin{proof}
    By Corollary~\ref{cor:1+1=2}, we have $\inv(D\rightarrow D)=2$ if and only $\inv(D)=1$, so the statement follows from Theorem \ref{thm:NP1}.
  \end{proof}

 \subsection{Solving  {\sc $k$-Tournament-Inversion} for $k\in \{1,2\}$}\label{subsec-poly}

 \begin{proposition}\label{prop:inv1}
 {\sc $1$-Tournament-Inversion} can be solved in $O(n^3)$ time.
 \end{proposition}
 \begin{proof}

 Let $T$ be a tournament.
  For every vertex $v$ one can check whether there is an inversion that transforms $T$ into a transitive tournament with source $v$.
 Indeed the unique possibility inversion is the one on the closed in-neighbourhood of $v$, $N^-[v]=N^-(v)\cup \{v\}$.
 So one can make inversion on $N^-[v]$ and check whether the resulting tournament is transitive.
 This can obviously be done in $O(n^2)$  time

 Doing this for every vertex $v$ yields an algorithm which solves {\sc $1$-Tournament-Inversion} in $O(n^3)$ time.
 \end{proof}

 \begin{theorem}\label{thm:inv2}
 {\sc $2$-Tournament-Inversion} can be solved in in $O(n^6)$ time.
 \end{theorem}

 The main idea to prove this theorem is to consider every pair $(s,t)$ of distinct vertices and to check whether there are two sets $X_1, X_2$ such that the inversion of $X_1$ and $X_2$ results in a transitive tournament with source $s$ and sink $t$.
 We need some definitions and lemmas.

 The {\bf symmetric difference} of two sets $A$ and $B$ is $A\triangle B= (A\setminus B) \cup (B\setminus A)$.

Let $T$ be a tournament and let $s$ and $t$ be two distinct vertices of $T$.
We define the following four sets
\begin{eqnarray*}
A(s,t) & = &  N^+(s) \cap N^-(t) \\
B(s,t) & = &  N^-(s) \cap N^+(t) \\
C(s,t) & = &  N^+(s) \cap N^+(t) \\
D(s,t) & = &  N^-(s) \cap N^-(t) \\
\end{eqnarray*}

 \begin{lemma}\label{lem:X1and2}
 Let $T$ be a tournament and let $s$ and $t$ be two distinct vertices of $T$.
 Assume there are two sets $X_1, X_2$ such that the inversion of $X_1$ and $X_2$ results in a transitive tournament with source $s$ and sink $t$.
 \begin{itemize}
\item[$(1)$] If $\{s,t\} \subseteq X_1\setminus X_2$, then $ts\in A(T)$, $C(s,t) = D(s,t) =\emptyset$ and $X_1=\{s,t\}\cup B(s,t)$.

 \item[$(2)$]  If $s\in X_1\setminus X_2$, $t\in X_2\setminus X_1$, then $st\in A(T)$, $A(s,t)\cap (X_1\cup X_2)=\emptyset$, $X_1=\{s\} \cup B(s,t) \cup D(s,t)$, and $X_2=\{t\}\cup B(s,t) \cup C(s,t)$.

 \item[$(3)$] If $s\in X_1\cap X_2$ and $t\in X_1\setminus X_2$, then $ts\in A(T)$, $X_1=\{s,t\}\cup B(s,t) \cup C(s,t)$, and $X_2=\{s\}\cup C(s,t) \cup D(s,t)$.

 \item[$(4)$] If $\{s,t\} \subseteq X_1\cap X_2$, then $st\in A(T)$, $C(s,t)=\emptyset$, $D(s,t)=\emptyset$, $X_1\cap X_2\subseteq A(s,t)\cup \{s,t\}$, and
 $B(s,t)=X_1\triangle X_2$.

 \end{itemize}
 \end{lemma}
 \begin{proof}
 (1) The arc between $s$ and $t$ is reversed once, so $ts\in A(T)$.

Assume for a contradiction, that there is a vertex $c\in C(S,t)$. The arc $tc$ must be reversed, so $c\in X_1$, but then the arc $sc$ is reversed contradicting the fact that $s$ becomes a source.
 Hence $C(s,t)=\emptyset$. Similarly $D(s,t)=\emptyset$. 
 
 The arcs from $t$ to $B(s,t)$ and from $B(s,t)$ to $s$ are reversed so $B(s,t)\subseteq X_1$.
 The arcs from $s$ to $A(s,t)$ and from $A(s,t)$ to $t$ are not reversed so $A(s,t) \cap X_1 =\emptyset$. 
 Therefore $X_1=\{s,t\}\cup B(s,t)$.
 
 \medskip
 
 (2) The arc between $s$ and $t$ is not reversed, so $st\in A(T)$.
 The arcs from $s$ to $A(s,t)$ and from $A(s,t)$ to $t$ are not reversed so $A(s,t) \cap X_1 =\emptyset$ and $A(s,t) \cap X_2 =\emptyset$. 
 The arcs from $t$ to $B(s,t)$ and from $B(s,t)$ to $s$ are reversed so $B(s,t)\subseteq X_1$ and $B(s,t)\subseteq X_2$.
 The arcs from $s$ to $C(s,t)$ are not reversed so $C(s,t)\cap X_1=\emptyset$ and the arcs from $t$ to $C(s,t)$ are reversed so $C(s,t)\subseteq X_2$.
 The arcs from $D(s,t)$ to $s$ are reversed so $D(s,t)\subseteq X_1$ and the arcs from $D(s,t)$ to $d$ are not reversed so $D(s,t)\cap X_2=\emptyset$.
 Consequently, $X_1=\{s\} \cup B(s,t) \cup D(s,t)$, and $X_2=\{t\}\cup B(s,t) \cup C(s,t)$.

 \medskip
 
 (3) The arc between $s$ and $t$ is reversed, so $ts\in A(T)$.
 The arcs from $A(s,t)$ to $t$ are not reversed so $A(s,t) \cap X_1 =\emptyset$. The arcs from $s$ to $A(s,t)$ are not reversed so $A(s,t) \cap X_2 =\emptyset$. 
  The arcs from $t$ to $B(s,t)$ are reversed so $B(s,t)\subseteq X_1$. The arcs from $B(s,t)$ to $s$ are reversed (only once) so $B(s,t) \cap X_2 =\emptyset$.
  The arcs from $t$ to $C(s,t)$ are reversed so $C(s,t)\subseteq X_1$. The arcs from $s$ to $C(s,t)$ must the be reversed twice so $C(s,t)\subseteq X_2$.
 The arcs from $D(s,t)$ to $t$ are not reversed so $D(s,t) \cap X_1 =\emptyset$.   The arcs from $D(s,t)$ to $s$ are reversed so $D(s,t)\subseteq X_2$.
  Consequently, $X_1=\{s,t\}\cup B(s,t) \cup C(s,t)$, and $X_2=\{s\}\cup C(s,t) \cup D(s,t)$.

\medskip 

(4) The arc between $s$ and $t$ is reversed twice, so $st\in A(T)$.

Assume for a contradiction, that there is a vertex $c\in C(s,t)$. The arc $tc$ must be reversed, so $c$ is in exactly one of $X_1$ ad $X_2$. But then the arc $sc$ is reversed contradicting the fact that $s$ becomes a source.
 Hence $C(s,t)=\emptyset$. Similarly $D(s,t)=\emptyset$. 
The arcs from $s$ to $A(s,t)$ and from $A(s,t)$ to $t$ are not reversed so every vertex of $A(s,t)$ is either in $X_1\cap X_2$ or in $V(T)\setminus (X_1\cup X_2)$. 
The arcs from $t$ to $B(s,t)$ and from $B(s,t)$ to $s$ are reversed so every vertex of $B(s,t)$ is either in $X_1\setminus X_2$ or in $X_2\setminus X_1$. 
Consequently, $X_1\cap X_2\subseteq A(s,t)\cup \{s,t\}$, and
 $B(s,t)=X_1\triangle X_2$.
 \end{proof}

  \begin{lemma}\label{lem:guess}
 Let $T$ be a tournament of order $n$ and let $s$ and $t$ be two distinct vertices of $T$.
 
 \begin{itemize}
\item[$(1)$] One can decide in $O(n^3)$ time whether there are two sets $X_1, X_2$ such that the inversion of $X_1$ and $X_2$ results in a transitive tournament with source $s$ and sink $t$ and $\{s,t\} \subseteq X_1\setminus X_2$. 

 \item[$(2)$]  One can decide in $O(n^2)$ time whether there are two sets $X_1, X_2$ such that the inversion of $X_1$ and $X_2$ results in a transitive tournament with source $s$ and sink $t$ and $s\in X_1\setminus X_2$ and $t\in X_2\setminus X_1$.

 \item[$(3)$] One can decide in $O(n^2)$ time whether there are two sets $X_1, X_2$ such that the inversion of $X_1$ and $X_2$ results in a transitive tournament with source $s$ and sink $t$ and $s\in X_1\cap X_2$ and $t\in X_1\setminus X_2$. 
 
  \item[$(4)$] One can decide in $O(n^4)$ time whether there are two sets $X_1, X_2$ such that the inversion of $X_1$ and $X_2$ results in a transitive tournament with source $s$ and sink $t$ and $\{s,t\} \subseteq X_1\cap X_2$. 
   \end{itemize}
 \end{lemma}
\begin{proof}
For all cases, we first compute $A(s,t)$, $B(s,t)$, $C(s,t)$, and $D(s,t)$, which can obviously be done in $O(n^2)$.

\medskip

(1) By Lemma~\ref{lem:X1and2}, we must have $ts\in A(T)$ and $C(s,t) = D(s,t) =\emptyset$.
So we first check if this holds.
Furthermore, by Lemma~\ref{lem:X1and2}, we must have $X_1=\{s,t\} \cup B(s,t)$. Therefore we invert  $\{s,t\} \cup B(s,t)$ which results in a tournament $T'$. Observe that $s$ is a source of $T'$ and $t$ is a sink of $T'$.
Hence, we return `Yes' if and only if $\inv(T'-\{s,t\})=1$ which can be tested in $O(n^3)$ by Proposition~\ref{prop:inv1}.

\medskip

(2) By Lemma~\ref{lem:X1and2}, we must have $st\in A(T)$. So we first check if this holds.
Furthermore,  by Lemma~\ref{lem:X1and2}, the only possibility is that $X_1=\{s\} \cup B(s,t) \cup D(s,t)$, and $X_2=\{t\}\cup B(s,t) \cup C(s,t)$.
So we invert those two sets and check whether the resulting tournament is a transitive tournament with source $s$ and sink $t$.
This can done in $O(n^2)$.

\medskip

(3) By Lemma~\ref{lem:X1and2}, we must have $ts\in A(T)$.
So we first check if this holds.
Furthermore, by Lemma~\ref{lem:X1and2}, the only possibility is that $X_1=\{s,t\}\cup B(s,t) \cup C(s,t)$, and $X_2=\{s\}\cup C(s,t) \cup D(s,t)$.
So we invert those two sets and check whether the resulting tournament is a transitive tournament with source $s$ and sink $t$.
This can done in $O(n^2)$.

\medskip

(4) By Lemma~\ref{lem:X1and2}, we must have $st\in A(T)$, $C(s,t)=\emptyset$, $D(s,t)=\emptyset$.
So we first check if this holds.
Furthermore, by Lemma~\ref{lem:X1and2}, the desired sets $X_1$ and $X_2$ must satisfy
$X_1\cap X_2\subseteq A(s,t)\cup \{s,t\}$, and
 $B(s,t)=X_1\triangle X_2$.
 
 In particular, every arc of $T_A=T\langle A(s,t)\rangle$ is either not reversed or reversed twice (which is the same). Hence $T_A$ must be a transitive tournament. So we check whether  $T_A$ is a transitive tournament and if yes, we find a directed hamiltonian path $P_A=(a_1, \dots , a_p)$ of it.
 This can be done in $O(n^2)$.
 
 Now we check that $B(s,t)$ admits a partition $(X'_1, X'_2)$ with $X'_i=X_i\cap B$ and the inversion of both $X'_1$ and $X'_2$ transforms $T\langle B(s,t)\rangle$ into a transitive tournament $T_B$ with source $s'$ and sink $t'$.
The idea is to investigate all possibilities for $s'$, $t'$ and the sets $X'_1$ and $X'_2$.
Since $(X'_1, X'_2)$ is a partition of $B(s,t)$ and $(X'_1,X'_2)$ is a decycling family if and only if $(X'_2,X'_1)$ is a decycling family, we may assume that

\begin{itemize}
\item[(a)] $\{s',t'\} \subseteq X'_1\setminus X'_2$, or 
\item[(b)]  $s'\in X'_1\setminus X'_2$ and $t'\in X'_2\setminus X'_1$.
\end{itemize}

For the possibilities corresponding to Case (a), we proceed as in (1) above.
For every arc $t's'\in A(T\langle B(s,t)\rangle)$, we check that $C(s',t') = D(s',t') =\emptyset$ (where those sets are computed in $T\langle B(s,t)\rangle$).
Furthermore, by Lemma~\ref{lem:X1and2}, we must have $X'_1=\{s,t\} \cup B(s',t')$ and $X'_2= B(s,t) \setminus X'_1$.
So we invert those two sets and check whether the resulting tournament $T_B$ is transitive. This can be done in $O(n^2)$ (for each arc $t's'$).

For the possibilities corresponding to Case (b), we proceed as in (2) above.
For every arc $t's'\in A(T\langle B(s,t)\rangle)$, by Lemma~\ref{lem:X1and2}, the only possibility is that $X'_1=\{s'\} \cup B(s',t') \cup D(s',t')$, and $X_2=\{t'\}\cup B(s',t') \cup C(s',t')$.
As those two sets form a partition of $B(s,t)$, we also must have $B(s',t')=\emptyset$ and $A(s',t')=\emptyset$.
So we invert those two sets and check whether the resulting tournament $T_B$ is transitive. This can be done in $O(n^2)$ for each arc $t's'$.

In both cases, we are left with a transitive tournament $T_B$.  We compute its directed hamiltonian path $P_B=(b_1, \dots , b_q)$ which can be done in $O(n^2)$.
We need to check whether this partial solution on $B(s,t)$ is compatible with the rest of the tournament, that is $\{s,t\} \cup A(s,t)$.
It is obvious that it will always be compatible with $s$ and $t$ as they become source and sink.
So we have to check that we can merge $T_A$ and $T_B$ into a transitive tournament on $A(s,t)$ and $B(s,t)$ after the reversals of $X_1$ and $X_2$.
In other words, we must interlace the vertices of $P_A$ and $P_B$. Recall that $Z=X_1\cap X_2\setminus \{s,t\}\subseteq A(s,t)$ and $X_i=X'_i\cup Z\cup \{s,t\}$, $i\in [2]$ so the arcs between $Z$ and $B(s,t)$ will be reversed exactly once when we invert $X_1$ and $X_2$. Using this fact,
one easily checks that this is possible if and only there are integers $j_1 \leq \dots \leq j_p$ such that
\begin{itemize}
\item either $b_j\ra a_i$ for $j\leq j_i$ and $b_j\la a_i$ for $j>j_i$ (in which case $a_i\notin Z$ and the arcs between $a_i$ and $B(s,t)$ are not reversed),
\item or  $b_j\la a_i$ for $j\leq j_i$ and $b_j\ra a_i$ for $j>j_i$ (in which case $a_i\in Z$ and the arcs between $a_i$ and $B(s,t)$ are reversed).

\end{itemize}
\noindent{} See Figure \ref{fig:merge} for an illustration of a case when we can merge the two orderings after reversing $X_1$ and $X_2$.

\begin{figure}[hbtp]
  \begin{center}
    \begin{tikzpicture}[scale=1]
      \tikzstyle{vertexB}=[circle,draw,  minimum size=10pt, scale=0.6, inner sep=0.5pt]
      \tikzstyle{vertexR}=[circle,draw,top color=red!60, bottom color=red!30,minimum size=10pt, scale=0.6, inner sep=0.5pt]
      \node (a1) at (1,4) [vertexB] {$a_1$};
      \node (a2) at (2,4) [vertexR] {$a_2$};
      \node (a3) at (3,4) [vertexB] {$a_3$};
      \node (a4) at (4,4) [vertexB] {$a_4$};
      \node (a5) at (5,4) [vertexB] {$a_5$};
      \node (a6) at (6,4) [vertexR] {$a_6$};
      \node (a7) at (7,4) [vertexB] {$a_7$};
      \node (a8) at (8,4) [vertexB] {$a_8$};
      \node (a9) at (9,4) [vertexB] {$a_9$};
      \node (a10) at (10,4) [vertexR] {$a_{10}$};
      \node (a11) at (11,4) [vertexB] {$a_{11}$};
      \draw[->,line width=0.03cm] (a1) to (a2);
      \draw[->,line width=0.03cm] (a2) to (a3);
      \draw[->,line width=0.03cm] (a3) to (a4);
      \draw[->,line width=0.03cm] (a4) to (a5);
      \draw[->,line width=0.03cm] (a5) to (a6);
      \draw[->,line width=0.03cm] (a6) to (a7);
      \draw[->,line width=0.03cm] (a7) to (a8);
      \draw[->,line width=0.03cm] (a8) to (a9);
      \draw[->,line width=0.03cm] (a9) to (a10);
      \draw[->,line width=0.03cm] (a10) to (a11);
      \draw (0.8,3.5) rectangle (4.2,4.5);
      \draw (4.8,3.5) rectangle (8.2,4.5);
      \draw (8.8,3.5) rectangle (11.2,4.5);
      
      \node (b1) at (1,0) [vertexB] {$b_1$};
      \node (b2) at (2,0) [vertexB] {$b_2$};
      \node (b3) at (3,0) [vertexB] {$b_3$};
      \node (b4) at (4,0) [vertexB] {$b_4$};
      \node (b5) at (5,0) [vertexB] {$b_5$};
      \node (b6) at (6,0) [vertexB] {$b_6$};
      \node (b7) at (7,0) [vertexB] {$b_7$};
      \node (b8) at (8,0) [vertexB] {$b_8$};
      \node (b9) at (9,0) [vertexB] {$b_9$};
      \node (b10) at (10,0) [vertexB] {$b_{10}$};
      \node (b11) at (11,0) [vertexB] {$b_{11}$};
      \node (b12) at (12,0) [vertexB] {$b_{12}$};
      \draw[->,line width=0.03cm] (b1) to (b2);
      \draw[->,line width=0.03cm] (b2) to (b3);
      \draw[->,line width=0.03cm] (b3) to (b4);
      \draw[->,line width=0.03cm] (b4) to (b5);
      \draw[->,line width=0.03cm] (b5) to (b6);
      \draw[->,line width=0.03cm] (b6) to (b7);
      \draw[->,line width=0.03cm] (b7) to (b8);
      \draw[->,line width=0.03cm] (b8) to (b9);
      \draw[->,line width=0.03cm] (b9) to (b10);
      \draw[->,line width=0.03cm] (b10) to (b11);
      \draw[->,line width=0.03cm] (b11) to (b12);
      \draw (0.8,-0.5) rectangle (3.2,0.5);
      \draw (3.8,-0.5) rectangle (6.2,0.5);
      \draw (6.8,-0.5) rectangle (9.2,0.5);
      \draw (9.8,-0.5) rectangle (12.2,0.5);
      \draw [->,line width=0.06cm, color=blue] (3.2,0.5) to (0.8,3.5);
      \draw [->,line width=0.06cm, color=blue] (4.2,3.5) to (3.8,0.5);
      \draw [->,line width=0.06cm, color=blue] (6.2,0.5) to (4.8,3.5);
      \draw [->,line width=0.06cm, color=blue] (8.2,3.5) to (6.8,0.5);
      \draw [->,line width=0.06cm, color=blue] (9.2,0.5) to (8.8,3.5);
      \draw [->,line width=0.06cm, color=blue] (11.2,3.5) to (9.8,0.5);
      \draw [->, line width=0.03cm, color=red] (a2) to (2,0.5);
      \draw [->, line width=0.03cm, color=red] (5,0.5) to (a2); 
      \draw [->, line width=0.03cm, color=red] (8,0.5) to (a2);
      \draw [->, line width=0.03cm, color=red] (11,0.5) to (a2);
      \draw [->, line width=0.03cm, color=red] (a6) to (2.1,0.5);
      \draw [->, line width=0.03cm, color=red] (a6) to (5.1,0.5); 
      \draw [->, line width=0.03cm, color=red] (8,0.5) to (a6);
      \draw [->, line width=0.03cm, color=red] (11,0.5) to (a6);
      \draw [->, line width=0.03cm, color=red] (a10) to (2.2,0.5);
      \draw [->, line width=0.03cm, color=red] (a10) to (5.2,0.5); 
      \draw [->, line width=0.03cm, color=red] (a10) to (8.2,0.5);
      \draw [->, line width=0.03cm, color=red] (11,0.5) to (a10);
      \draw [->, line width=0.06cm] (5.5,5) to (7.5,5);
      \draw [->, line width=0.06cm] (5,-1) to (8,-1);
      
      \end{tikzpicture}
  \end{center}
  \caption{Indicating how to merge the two orderings of $A$ and $B$. The fat blue edges indicate that the final ordering will be $b_1-b_3,a_1-a_4,b_4-b_6,a_5-a_8,b_7-b_9,a_9-a_{11},b_{10}-b_{12}$. The set $Z=\{a_2,a_6,a_{10}\}$ consists of those  vertices from $A(s,t)$ which are in $X_1\cap X_2$. These vertices are shown in red. The red arcs between a vertex of $Z$  and one of the boxes indicate that all arcs between the vertex and those of the box have the direction shown. Hence the boxes indicate that values of  $j_1,\ldots{},j_{11}$ satisfy that : $j_1=\ldots{}=j_4=3$, $j_5=\ldots{}=j_8=6$, $j_9=\ldots{}=j_{11}=9$.}\label{fig:merge}
  \end{figure}
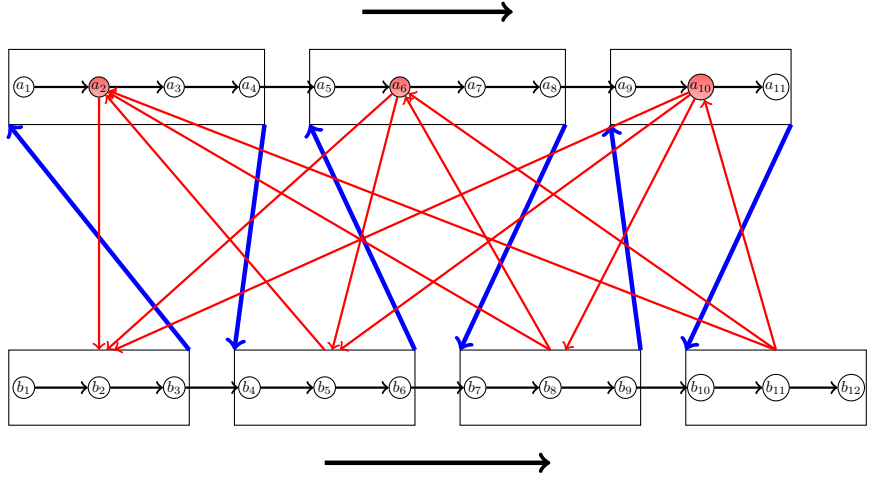

Deciding whether there are such indices can be done in $O(n^2)$ for each possibility.

As we have $O(n^2)$ possibilities, and for each possibility the procedure runs in $O(n^2)$ time, the overall procedure runs in $O(n^4)$ time.
\end{proof}

 \begin{proof}[Proof of Theorem~\ref{thm:inv2}]
 By Lemma~\ref{lem:reduc}, by removing iteratively the sources and sinks of the tournament, it suffices to solve the problem for a tournament with  no sink and no source.
 
Now for each pair $(s,t)$ of distinct vertices, one shall check whether there are two sets $X_1, X_2$ such that the inversion of $X_1$ and $X_2$ results in a transitive tournament with source $s$ and sink $t$.
Observe that since $s$ and $t$ are neither sources nor sinks in $T$, each of them must belong to at least one of $X_1,X_2$.
Therefore, without loss of generality, we are in one of the following possibilities:
\begin{itemize}
\item $\{s,t\} \subseteq X_1\setminus X_2$. Such a possibility can be checked in $O(n^3)$ by Lemma~\ref{lem:guess} (1).
\item $s\in X_1\setminus X_2$ and $t\in X_2\setminus X_1$. Such a possibility can be checked in $O(n^2)$ by Lemma~\ref{lem:guess} (2).
\item $s\in X_1\cap X_2$ and $t\in X_1\setminus X_2$. Such a possibility can be checked in $O(n^2)$ by Lemma~\ref{lem:guess} (3).
\item $t\in X_1\cap X_2$ and $s\in X_1\setminus X_2$. Such a possibility is the directional dual of the preceding one. It can be tested  in $O(n^2)$ by reversing all arcs and applying Lemma~\ref{lem:guess} (3).
\item  $\{s,t\} \subseteq X_1\cap X_2$. Such a possibility can be checked in $O(n^4)$ by Lemma~\ref{lem:guess} (4).
\end{itemize}

Since there are $O(n^2)$ pairs $(s,t)$ and for each pair the procedure runs in $O(n^4)$, the algorithm runs in $O(n^6)$ time.
 \end{proof}

 \subsection{Computing related parameters when the inversion number is bounded}\label{subsec:bounded-inv}

 The aim of this subsection is to prove the following theorem.
 
 \begin{theorem}\label{thm:NP-inv1}
 Let $\gamma$ be a parameter in $\tau, \tau', \nu$.
Given an oriented graph $D$ with inversion number $1$ and an integer $k$, it is NP-complete to decide whether $\gamma(D) \leq k$.
\end{theorem}

 Let $D$ be a digraph.
The {\bf second subdivision} of $D$ is the oriented graph $S_2(D)$ obtained from $D$ by replacing every arc $a=uv$ by a directed path $P_a=(u, x_a, y_a, u)$ where $x_a, y_a$ are two new vertices. 

\begin{lemma}\label{lem:subdiv}
Let $D$ be a digraph.
\begin{itemize}
\item[(i)] $\inv(S_2(D)) \leq 1$.
\item[(ii)] $\tau'(S_2(D)) = \tau'(D)$,  $\tau(S_2(D)) = \tau(D)$, and $\nu(S_2(D)) = \nu(D)$.
\end{itemize}
\end{lemma}
\begin{proof}
(i) Inverting the set $\bigcup_{a\in A(D)} \{x_a, y_a\}$ makes $S_2(D)$ acyclic. Indeed the $x_a$ become sinks, the $y_a$ become sources and
the other vertices form a stable set. Thus $\inv(S_2(D))=1$.

\medskip

(ii) There is a one-to-one correspondence between directed cycles in $D$ and directed cycles in $S_2(D)$ (their second subdivision).
Hence $\nu(S_2(D)) = \nu(D)$.

Moreover every cycle transversal $S$ of $D$ is also a cycle transversal of $S_2(D)$. So  $\tau(S_2(D)) \leq \tau(D)$.
Now consider a  cycle transversal $T$. If  $x_a$ or $y_a$ is in $S$ for some $a\in A(D)$, then we can replace it by any end-vertex of $a$. 
Therefore, we may assume that $T\subseteq V(D)$, and so $T$ is a cycle transversal of $D$.
Hence $\tau(S_2(D)) = \tau(D)$.

Similarly, consider a cycle arc-transversal $F$ of $D$. Then $F'= \{a \mid x_ay_a\in F\}$ is a cycle arc-transversal of $S_2(D)$.
Conversely, consider a cycle arc-transversal  $F'$ of $S_2(D)$. Replacing each arc incident to $x_a$, $y_a$ by $x_ay_a$ for each $a\in A(D)$, we obtain another cycle arc-transversal. So we may assume that $F'\subseteq \{x_ay_a \mid a\in A(D)\}$. Then $F= \{a \mid x_ay_a\in F'\}$ is a cycle arc-transversal of $D$.
Thus  $\tau'(S_2(D)) = \tau'(D)$.
\end{proof}

  \begin{proof}[Proof of Theorem~\ref{thm:NP-inv1}]
Since computing each of $\tau$, $\tau'$, $\nu$ is NP-hard, Lemma~\ref{lem:subdiv}~(ii) implies that computing each of $\tau$, $\tau'$, $\nu$ is also NP-hard for second subdivisions of digraphs. As those oriented graphs have inversion number $1$ (Lemma~\ref{lem:subdiv}~(i)), computing each of $\tau$, $\tau'$, $\nu$ is NP-hard for oriented graphs with inversion number $1$.
\end{proof}

\acknowledgements
The authors would like to thank Maurice Pouzet for introducing the topic of inversion to us. They are very pleased to honour him on his 75th birthday.
They are also grateful to an anonymous referee whose remarks helped to improve the presentation of the paper.

This research was supported by the Independent Research Fund Denmark under grant number DFF 7014-00037B, and by the french Agence Nationale de la Recherche under contract Digraphs ANR-19-CE48-0013-01.
This work was done while J\o{}rgen Bang-Jensen was visiting  team Coati, I3S and INRIA Sophia Antipolis. Hospitality and financial support is gratefully acknowledged. Ce travail a b\'en\'efici\'e d'une aide du gouvernement fran\c{c}ais, 
g\'er\'ee par l'Agence Nationale de la Recherche au titre du projet Investissements d'Avenir UCAJEDI portant la r\'ef\'erence no ANR-15-IDEX-01.
This work was done while Jonas Costa Ferreira da Silva was spending his second year of PhD as ``sandwich'' in team Coati, I3S and INRIA Sophia Antipolis and this visit was funded by the project CAPES STIC-AmSud GALOP 88881.197438/2018-01.


\end{document}